\documentclass[11pt,a4paper]{amsart}
\usepackage[T1]{fontenc}
\usepackage{lmodern}
\usepackage[margin=30mm]{geometry}
\usepackage{amsmath,amssymb,mathtools}
\usepackage{esint}
\usepackage{enumitem}
\usepackage{microtype}
\usepackage{xcolor}
\usepackage{url}
\usepackage[colorlinks=true,linkcolor=blue!45!black,citecolor=blue!45!black,urlcolor=blue!45!black]{hyperref}
\hypersetup{pdftitle={Moser--Trudinger inequalities and the complex Monge--Amp\`ere mean field equation},pdfauthor={}}
\numberwithin{equation}{section}
\newtheorem{theorem}{Theorem}[section]
\newtheorem{proposition}[theorem]{Proposition}
\newtheorem{lemma}[theorem]{Lemma}
\newtheorem{corollary}[theorem]{Corollary}
\newtheorem*{theorem*}{Theorem}
\newtheorem*{conjecture}{Conjecture}
\theoremstyle{remark}
\newtheorem{remark}[theorem]{Remark}
\newtheorem*{remark*}{Remark}
\newcommand{\C}{\mathbb C}
\newcommand{\R}{\mathbb R}
\newcommand{\PSH}{\operatorname{PSH}}
\newcommand{\Ric}{\operatorname{Ric}}
\newcommand{\Capa}{\operatorname{Cap}}
\newcommand{\ddc}{dd^c}
\newcommand{\dV}{\,dV}
\newcommand{\dsigma}{\,d\sigma}
\newcommand{\pbar}{\bar\partial}

\setlist[enumerate]{label=\textup{(\roman*)},leftmargin=2.3em,itemsep=2pt}
\allowdisplaybreaks[2]
\title[Moser--Trudinger inequalities]{Moser--Trudinger inequalities and the complex Monge--Amp\`ere
mean field equation}

\author{Jie Deng}
\author{Jiaxiang Wang}
\author{Xu-jia Wang}
\author{Bin Zhou}

\address{School of Mathematical Sciences, Peking University, Beijing 100871, China}
\email{dj0401@stu.pku.edu.cn}
\address{School of Mathematical Sciences and LPMC, Nankai University, Tianjin 300071, China}
\email{wangjx\_manifold@126.com}
\address{Institute for Theoretical Sciences, Westlake University, Hangzhou, 310030, China}
\email{Wangxujia@westlake.edu.cn}
\address{School of Mathematical Sciences, Peking University, Beijing 100871, China}
\email{bzhou@pku.edu.cn}

\thanks{Jiaxiang Wang was partially supported by the National Natural Science Foundation of China (Grant
Nos. 12401068 and 12571060). Jie Deng and Bin Zhou were partially supported by the National Key R\&D Program of China (Grant No.~2023YFA009900) and the National Natural Science Foundation of China (Grant No.~12271008).}
\date{}
\begin{document}
\begin{abstract}
We prove the Moser--Trudinger inequality in exponential form at the critical coefficient for the complex Monge--Amp\`ere energy on smooth bounded strictly pseudoconvex domains. We also determine the optimal coefficients in the logarithmic Moser--Trudinger inequality on the unit ball, and obtain a volume--capacity inequality with the optimal coefficient in the exponent. The logarithmic inequality follows from radial symmetry and uniqueness for the associated mean field equation. We prove that every smooth solution has constant boundary normal derivative, giving existence and uniqueness precisely for total masses $0<\Lambda<[2\pi(n+1)]^n$ and proving Berman--Berndtsson's conjecture. We also establish smoothness of bounded pluripotential and finite-energy solutions. The volume--capacity estimate, a comparison of capacities of sublevel sets at different levels, and Moser's lemma give the critical exponential inequality on the ball. Domain monotonicity then transfers the inequality to general domains.
\end{abstract}
\maketitle

\section{Introduction}\label{sec:introduction}

Sobolev, Moser--Trudinger and related functional inequalities are fundamental tools in analysis, partial differential equations and geometry. They control the integrability of functions in terms of their energies and play a central role in the existence and regularity theory of nonlinear equations and in geometric variational problems.

The Moser--Trudinger inequality describes the exponential integrability that replaces a Sobolev embedding at the critical exponent. Let $\Omega\subset\mathbb R^n$, $n\geq2$, be a bounded domain. Trudinger \cite{Tr67} established the limiting exponential embedding, and Moser \cite{M71} proved its sharp form in every dimension $n\geq2$:
\begin{equation}\label{classical-MT}
\int_\Omega\exp\left[a_n
\left(\frac{|u|}{\|\nabla u\|_{L^n(\Omega)}}\right)^{\frac{n}{n-1}}\right]dx
\leq C_n|\Omega|
\end{equation}
for every nonzero $u\in W^{1,n}_0(\Omega)$, where
$a_n=n\omega_{n-1}^{\frac{1}{n-1}}$ and $\omega_{n-1}=|\mathbb S^{n-1}|$.
This is the critical exponential inequality associated with the $n$-Laplace operator
\[
\Delta_n u=\operatorname{div}\bigl(|\nabla u|^{n-2}\nabla u\bigr).
\]
The coefficient $a_n$ is optimal.

\vskip 10pt

A fully nonlinear analogue is furnished by the real $k$-Hessian operator
\[S_k(D^2u)=\sigma_k(\lambda(D^2u)),\]
where $\sigma_k$ is the $k$th elementary symmetric polynomial. Suppose that $\Omega\subset\mathbb R^n$ is a bounded smooth strictly $(k-1)$-convex domain. Denote by $\Phi_0^k(\Omega)$ the class of $k$-admissible functions vanishing on $\partial\Omega$; for smooth functions, admissibility means $S_j(D^2u)\geq0$ for $1\leq j\leq k$. Define the Hessian energy and Hessian norm by
\[
\mathcal I_k(u)=\frac{1}{k+1}\int_\Omega(-u)S_k(D^2u)\,dx,
\qquad \|u\|_{\Phi_0^k(\Omega)}=\mathcal I_k(u)^{\frac{1}{k+1}}.
\]
Wang \cite{W94} established Sobolev inequalities for the $k$-Hessian operator. For $n=2k$,
Tian and Wang \cite{TiW10} later proved the sharp Moser--Trudinger inequality for every nonzero $u\in\Phi_0^k(\Omega)$:
\begin{equation}\label{hessian-MT}
\int_\Omega\exp\left[b_{n,k}
\left(\frac{-u}{\|u\|_{\Phi_0^k(\Omega)}}\right)^{\frac{k+1}{k}}\right]dx
\leq C,
\end{equation}
where $C=C(n,\Omega)$, $k=\frac{n}{2}$ and
\[b_{n,k}=n\left[\frac{\omega_{n-1}}{k(k+1)}\binom{n-1}{k-1}\right]^{\frac{1}{k}}=n\left[\frac{\omega_{n-1}}{n(k+1)}\binom{n}{k}\right]^{\frac{2}{n}}.\]
For $n=2$ and $k=1$, \eqref{hessian-MT} recovers the planar inequality \eqref{classical-MT} on the subharmonic class.

\vskip 10pt
In several complex variables, the complex Monge--Amp\`ere operator provides a corresponding energy for plurisubharmonic functions. Let $\Omega\subset\mathbb C^n$ be a bounded strictly pseudoconvex domain with smooth boundary, and set
\[
\mathcal H_0(\Omega)
=\{u\in\PSH(\Omega)\cap C^\infty(\overline\Omega):
 u|_{\partial\Omega}=0\}.
\]
With the convention $\ddc=i\partial\bar\partial$, define the Monge--Amp\`ere energy and its associated norm by
\[
\mathcal E(u)=\frac{1}{n+1}\int_\Omega(-u)(\ddc u)^n,
\qquad
\|u\|_{\mathcal H_0(\Omega)}=\mathcal E(u)^{\frac{1}{n+1}}.
\]
The energy is homogeneous of degree $n+1$. For locally bounded plurisubharmonic functions, the measure $(\ddc u)^n$ is defined by Bedford--Taylor \cite{BT76,BT82}.

Moser--Trudinger inequalities associated with complex Monge--Amp\`ere energies have their origins in Aubin's work on K\"ahler--Einstein metrics \cite{Au84}. The domain setting has subsequently been studied through pluripotential theory and variational methods \cite{BB22,C19,GKY13,AC19}. The corresponding Moser--Trudinger inequality takes the form
\begin{equation}\label{complex-MT}
\int_\Omega
\exp\left[\alpha
\left(\frac{-u}{\|u\|_{\mathcal H_0(\Omega)}}\right)^{\beta}\right]\dV
\leq C,
\qquad 0\ne u\in\mathcal H_0(\Omega),
\end{equation}
where $dV$ denotes Euclidean volume. We refer to \eqref{complex-MT} as the
\emph{Moser--Trudinger inequality in exponential form}.
Berman and Berndtsson \cite{BB22} established quasi-sharp inequalities on pseudoconvex domains and related them to Brezis--Merle type estimates. Set
\[
 \alpha_n=\frac{2\pi n}{(n+1)^{\frac{1}{n}}}.
\]
Their result \cite[Theorem~1.5]{BB22} gives \eqref{complex-MT} with $\alpha=(1-\varepsilon)\alpha_n$ and $\beta=\frac{n+1}{n}$ for every $\varepsilon\in(0,1)$. More quantitatively, their logarithmic estimate in \cite[Section~4]{BB22} takes the following form in our normalization:
\[
 \log\int_\Omega e^{-v}\dV
 \leq\frac{\mathcal E(v)}{[2\pi(n+1)]^n(1-\eta)^n}
 +(n-1)\log\frac{1}{\eta}+\log K_{n,\Omega},
 \qquad 0<\eta<\frac{1}{2},
\]
where $K_{n,\Omega}$ is independent of $\eta$ and $v\in\mathcal H_0(\Omega)$.
This is a \emph{logarithmic Moser--Trudinger inequality}.
The threshold $\alpha_n$ is sharp: \eqref{complex-MT} fails for $\beta=\frac{n+1}{n}$ and $\alpha>\alpha_n$. In this paper we prove that the exponential inequality holds at $\alpha=\alpha_n$. We also determine both optimal constants in the logarithmic inequality on the ball and obtain a volume--capacity inequality with the optimal coefficient.

\subsection{Sharp Moser--Trudinger inequalities}\label{subsec:sharp-inequalities}
The rearrangement argument used in the classical theory does not extend directly to arbitrary plurisubharmonic functions: Schwarz symmetrization need not preserve plurisubharmonicity. Berman and Berndtsson \cite{BB14} proved that symmetrization decreases the Monge--Amp\`ere energy for $S^1$-invariant plurisubharmonic functions on the ball, and obtained the sharp Moser--Trudinger inequality in that setting.

A complementary approach was developed in our earlier work \cite{WWZ20} by classical PDE methods inspired by \cite{W94,TiW10}, which deal with real Hessian equations. The proof combines a descending flow for a Monge--Amp\`ere functional with Sobolev and generalized Brezis--Merle estimates. This establishes exponential integrability with the optimal power $\frac{n+1}{n}$, but does not determine the optimal coefficient in the exponent.

Our first result proves the exponential inequality at the critical coefficient, without a symmetry assumption.
\begin{theorem}\label{thm:sharp-MT}
Let $\Omega\Subset B_R(a)\subset\mathbb C^n$ be a strictly pseudoconvex domain with smooth boundary. For every nonzero $u\in\mathcal H_0(\Omega)$,
\[
\int_\Omega
\exp\!\left[
\alpha_n
\left(\frac{-u}{\|u\|_{\mathcal H_0(\Omega)}}\right)^{\frac{n+1}{n}}
\right]\dV
\leq C_n|B_R|\cdot\exp\!\left(\sum_{j=1}^n\frac{1}{j}\right),
\]
where $C_n$ depends only on $n$. The coefficient $\alpha_n$ is optimal.
\end{theorem}

We next state the sharp logarithmic inequality on the unit ball $B=\{z\in\C^n:|z|<1\}$.
We use Cegrell's classes \cite{C98,C04}. The class $\mathcal E_0(B)$ consists of bounded negative plurisubharmonic functions with zero boundary limits and finite Monge--Amp\`ere mass. A negative plurisubharmonic function $u$ belongs to $\mathcal E^1(B)$ if there exist $u_j\in\mathcal E_0(B)$ decreasing to $u$ such that
\[
 \sup_j\int_B(-u_j)(\ddc u_j)^n<\infty.
\]
The Monge--Amp\`ere measure and the finite energy $\mathcal E(u)$ are defined on $\mathcal E^1(B)$ by decreasing approximation with such $u_j$.

\begin{theorem}\label{thm:sharp-log}
Let $n\geq1$. For every $u\in\mathcal E^1(B)$,
\begin{equation}\label{eq:sharp-mt}
 \log\left(\frac{1}{|B|}\int_B e^{-u}\dV\right)
 \leq\frac{\mathcal E(u)}{[2\pi(n+1)]^n}
       +\sum_{j=1}^n\frac{1}{j}.
\end{equation}
In particular, this holds for every $u\in\mathcal H_0(B)$, without any symmetry assumption. The coefficient of $\mathcal E(u)$ is sharp. For this coefficient, the additive constant $\sum_{j=1}^n\frac{1}{j}$ is also optimal.
\end{theorem}

\subsection{Volume--capacity inequalities}\label{subsec:volume-capacity}
The logarithmic inequality also gives a sharp estimate for the volume of a set in terms of its relative Monge--Amp\`ere capacity. Such volume--capacity inequalities are also called \emph{isocapacitary inequalities}. For a compact set $K\Subset\Omega$, define
\[
 \Capa(K,\Omega)
 =\sup\left\{\int_K(\ddc v)^n:
 v\in\PSH(\Omega),\ -1\leq v\leq0\right\}.
\]
The volume--capacity inequality of {\AA}hag, Cegrell, Ko{\l}odziej, Ph{\d{a}}m and Zeriahi \cite{ACKPZ09} contains a polynomial factor in front of the exponential term. Removing this factor while retaining the optimal coefficient in the exponent was asked in \cite[Section~5.2]{DGZ16}. Our next result gives such an estimate on the ball.

\begin{theorem}\label{thm:volume-capacity}
For every compact set $K\Subset B$,
\begin{equation}\label{eq:sharp-volume-capacity}
 \frac{|K|}{|B|}
 \leq\exp\left(\sum_{j=1}^n\frac{1}{j}
 -\frac{2\pi n}{\Capa(K,B)^{\frac{1}{n}}}\right).
\end{equation}
The right-hand side is understood as zero when $\Capa(K,B)=0$.
The coefficient $2\pi n$ is optimal.
\end{theorem}

By monotonicity of relative capacity, Theorem~\ref{thm:volume-capacity} also gives a volume--capacity estimate on every bounded hyperconvex domain contained in a ball. The comparison principle then gives the following estimate for sublevel sets, which removes the polynomial factor in the uniform volume estimate of \cite{ACKPZ09} and answers the corresponding question in \cite[Section~5.2]{DGZ16}.

For a bounded hyperconvex domain $\Omega$, the class $\mathcal E_0(\Omega)$ is defined as above. Cegrell's class $\mathcal F(\Omega)$ consists of negative plurisubharmonic functions $u$ admitting a decreasing sequence $u_j\in\mathcal E_0(\Omega)$ with
\[
 \sup_j\int_\Omega(\ddc u_j)^n<\infty.
\]
\begin{corollary}\label{cor:sublevel-volume}
Let $\Omega\Subset B_R(a)\subset\mathbb C^n$ be a bounded hyperconvex domain. Then, for every nonzero $u\in\mathcal F(\Omega)$ and every $s>0$,
\begin{equation}\label{eq:sublevel-volume}
 \bigl|\{z\in\Omega:u(z)<-s\}\bigr|
 \leq |B_R|\cdot\exp\left(
 \sum_{j=1}^n\frac{1}{j}
 -\frac{2\pi n s}
 {\left(\displaystyle\int_\Omega(\ddc u)^n\right)^{\frac{1}{n}}}
 \right).
\end{equation}
The coefficient $2\pi n$ is optimal.
\end{corollary}

\subsection{The mean field equation and proof strategy}\label{subsec:mean-field}
We now introduce the mean field equation used to prove the logarithmic inequality. For $\beta=\frac{n+1}{n}$, Young's inequality and the homogeneity of $\mathcal E$ give the logarithmic form associated with \eqref{complex-MT}:
\begin{equation}\label{eq:log-energy-mt}
 \log\left(\frac{1}{|\Omega|}\int_\Omega e^{-u}\dV\right)
 \leq\frac{n^n}{(n+1)^{n+1}\alpha^n}\mathcal E(u)+C_\alpha.
\end{equation}
At $\alpha=\alpha_n$, the energy coefficient is $\frac{1}{[2\pi(n+1)]^n}$. Theorem~\ref{thm:sharp-log} determines the optimal additive constant on the ball for this coefficient.

On $\mathcal E^1(B)$, the logarithmic inequality leads to the functional
\[
 \mathcal G_\Lambda(u)
 =-\mathcal E(u)+\Lambda\log\int_B e^{-u}\dV,
 \qquad u\in\mathcal E^1(B),\quad\Lambda>0,
\]
whose Euler--Lagrange equation is the complex Monge--Amp\`ere mean field equation
\begin{equation}\label{eq:main}
 \begin{cases}
 (\ddc u)^n=\Lambda\frac{e^{-u}}{\int_B e^{-u}\dV}\,dV
       &\text{in }B,\\[5pt]
 u=0&\text{on }\partial B.
 \end{cases}
\end{equation}
Here $\Lambda=\int_B(\ddc u)^n$ is the total Monge--Amp\`ere mass. A bounded weak solution of \eqref{eq:main} is a function $u\in\PSH(B)\cap L^\infty(B)$ satisfying the equation in the Bedford--Taylor sense and
\[
 \lim_{B\ni z\to\zeta}u(z)=0,\qquad \zeta\in\partial B.
\]

The variational approach to such equations was developed in \cite{BBGZ13,B13}. The coefficient in \eqref{eq:log-energy-mt} identifies $[2\pi(n+1)]^n$ as the critical mass. In the subcritical range, the quasi-sharp inequality of Berman and Berndtsson already shows that $\mathcal G_\Lambda$ is bounded above and tends to $-\infty$ as $\mathcal E(u)\to\infty$. Their variational existence theorem \cite[Theorem~7.2]{BB22} gives a maximizer in $\mathcal E^1(B)$ satisfying \eqref{eq:main}.

Berman and Berndtsson made the following conjecture.
\begin{conjecture}[Berman--Berndtsson \cite{BB22}]
Every solution of \eqref{eq:main} on the unit ball is radial.
\end{conjecture}
The radial solutions have been classified in \cite{L16}, so the conjecture also gives uniqueness for each fixed $\Lambda$. Our normalization is related to that of \cite{BB22} by $u=\gamma\varphi$ and $\Lambda=(2\pi\gamma)^n$, with subcritical range $0<\gamma<n+1$.

In complex dimension one, the moving-plane theorem of Gidas, Ni and Nirenberg \cite{GNN79} gives radial symmetry on the unit disc. In higher complex dimension, a reflection in a real hyperplane generally preserves neither plurisubharmonicity nor the complex Monge--Amp\`ere operator, so the moving-plane method does not apply directly.
Delano\"e \cite{D85} studied radial boundary value problems for real and complex Monge--Amp\`ere equations. Guedj, Kolev and Yeganefar \cite{GKY13} proved uniqueness among $S^1$-invariant positive K\"ahler--Einstein fillings when the domain is strictly convex with respect to the solution metric. Lu and Phung \cite[Theorem~1.2]{LP25} proved uniqueness of bounded pluripotential solutions on bounded hyperconvex domains, without symmetry assumptions, for sufficiently small positive temperature and $L^p$ densities with $p>1$.

We prove that every smooth solution has constant boundary normal derivative. This gives radial symmetry, and the known radial classification then gives uniqueness. We also treat bounded weak solutions and finite-energy solutions.
\begin{theorem}\label{thm:MT}
Let $n\geq1$. Equation \eqref{eq:main} has a smooth plurisubharmonic solution on $\overline B$ if and only if
\begin{equation}\label{eq:range}
 0<\Lambda<[2\pi(n+1)]^n.
\end{equation}
Moreover, the solution is unique and given by
\begin{equation}\label{eq:explicit}
 u_\Lambda(z)=(n+1)\log\frac{1+\mu|z|^2}{1+\mu},
 \qquad
 \mu=\frac{\Lambda^{\frac{1}{n}}}{2\pi(n+1)-\Lambda^{\frac{1}{n}}}.
\end{equation}
\end{theorem}
Section~\ref{sec:smoothness} shows that every bounded weak solution is smooth on $\overline B$ and that every finite-energy solution is bounded. Consequently, Theorem~\ref{thm:MT}, including its nonexistence assertion, holds for both of these weak solution classes.

The proofs of the inequalities proceed in the following order. Theorem~\ref{thm:MT} identifies the subcritical maximizers of $\mathcal G_\Lambda$ with the explicit radial family. Passing to the critical mass gives Theorem~\ref{thm:sharp-log}, including the optimal additive constant. Applying this logarithmic inequality to multiples of the relative extremal function gives Theorem~\ref{thm:volume-capacity}. A comparison of capacities of sublevel sets at different levels then gives a one-dimensional energy estimate, to which we apply Moser's lemma to obtain the critical exponential inequality on the ball. Finally, domain monotonicity gives Theorem~\ref{thm:sharp-MT} on general domains.

\smallskip
\noindent\textbf{Outline of the paper.}
Section~\ref{sec:preliminaries} gives the basic identities on the ball and recalls the classification of radial solutions. Section~\ref{sec:smoothness} proves smoothness of bounded weak solutions and finite-energy solutions. Section~\ref{sec:normal} proves constancy of the boundary normal derivative and completes the proof of Theorem~\ref{thm:MT}. In Section~\ref{sec:sharp-log-ball} we prove Theorem~\ref{thm:sharp-log}. Section~\ref{sec:exponential-ball} proves the volume--capacity and sublevel estimates and derives the critical exponential inequality on the ball. Section~\ref{sec:general-domains} proves domain monotonicity and completes the proof on general domains.

\section*{Acknowledgements and AI disclosure}

The preprint of Chen and Ma
\cite[Theorem~1.1 and Corollaries~1.2--1.3]{CM26}, which appeared
online on November~17, 2026, proves radial symmetry
and classifies smooth solutions of the complex Monge--Amp\`ere mean
field equation on the ball. It contains a similar proof of the
sharp logarithmic Moser--Trudinger inequality. Our work also establishes
smoothness of bounded weak solutions and finite-energy solutions on
the ball, and proves domain monotonicity.

The proofs of Theorem~\ref{thm:MT} and
Theorem~\ref{thm:sharp-log} were obtained by the authors of the first
version by late September 2026. Bin Zhou presented these results
at the workshop on recent advances in Monge--Amp\`ere equations held
at Westlake University in Hangzhou on November~11, 2026.
Work on the Moser--Trudinger inequality in exponential form at the
critical coefficient was still ongoing when the first version was
posted. The present revision establishes this inequality and adds
a volume--capacity inequality with the optimal coefficient in the exponent.

During the preparation and revision of this paper, we used ChatGPT
(OpenAI) to explore auxiliary functions involving the gradient of the
solution. The function it suggested plays a key role in the proofs
of Lemma~\ref{lem:F-inequality} and Proposition~\ref{prop:normal}.
We verified the computations in Section~\ref{sec:preliminaries}
and carried out those in Section~\ref{sec:normal}.
ChatGPT was also used to develop and revise the capacity arguments in
Section~\ref{sec:exponential-ball}. Further uses included literature
searches, improvements to the exposition, and checks of formatting
and typographical errors. We take full responsibility for all
mathematical statements, proofs, computations and references in this paper.
\section{Notation and preliminaries}\label{sec:preliminaries}

\subsection{Notation and basic identities}\label{subsec:notation}
We use
\[
 d^c=\frac{i}{2}(\pbar-\partial),\qquad
 \ddc=i\partial\pbar,
 \qquad |B|=\frac{\pi^n}{n!},\qquad |\partial B|=2n|B|.
\]
For smooth $u$,
\[
 (\ddc u)^n=2^n n!\det(u_{j\bar k})\,dV.
\]
For a smooth strictly plurisubharmonic solution, set $g=(g_{j\bar k}):=(u_{j\bar k})$.

Let $(g^{j\bar k})$ denote the transpose of its inverse.
The cofactor tensor of the complex Hessian is divergence-free, i.e.,
\begin{equation}\label{eq:cofactor}
 \sum_j\partial_j\bigl(\det g\,g^{j\bar k}\bigr)=0.
\end{equation}
With repeated indices summed, define
\[
 \Delta_g f=g^{j\bar k}f_{j\bar k},
 \qquad |\partial f|_g^2=g^{j\bar k}f_jf_{\bar k}
\]
for $f\in C^\infty(\overline B)$.
In particular, $\Delta_g u=n$.

Let $\nabla_g$ be the Levi--Civita connection of the K\"ahler metric $g$.
We denote the $(2,0)$-part of the covariant Hessian of $u$ by
\[
 \mathcal H=\nabla_g^{1,0}\partial u
 =\mathcal H_{ij}\,dz^i\otimes dz^j.
\]
In holomorphic coordinates, $\mathcal H_{ij}$ is the $j$th component
of $\nabla_i(\partial u)$, the covariant derivative of $\partial u$
in the $z^i$ direction. Thus
\[
 \mathcal H_{ij}=\nabla_i u_j
 =u_{ij}-g^{k\bar\ell}(\partial_i g_{j\bar\ell})u_k.
\]
We write $\mathcal H(\partial u)$ for the $(1,0)$-form obtained by
inserting the vector $g^{j\bar k}u_{\bar k}\frac{\partial}{\partial z^j}$
in the second slot of $\mathcal H$. Its components are
\[
 [\mathcal H(\partial u)]_i=\mathcal H_{ij}g^{j\bar k}u_{\bar k}.
\]

When $u$ solves \eqref{eq:main}, write $Z=\displaystyle\int_B e^{-u}\dV$. Then
\begin{equation}\label{eq:det-equation}
 \det g=\frac{\Lambda}{2^n n!Z}e^{-u}.
\end{equation}

\begin{lemma}[Bochner formula]\label{lem:bochner}
Let $u$ be a smooth strictly plurisubharmonic solution of \eqref{eq:main}.
Then
\[
 \Delta_g\bigl(|\partial u|_g^2\bigr)=n+|\mathcal H|_g^2+|\partial u|_g^2.
\]
\end{lemma}
\begin{proof}
Equation~\eqref{eq:det-equation} implies that $\log\det g+u$ is constant, so
\[
 \Ric(g)_{i\bar j}=-\partial_i\partial_{\bar j}\log\det g=g_{i\bar j}.
\]
Since $\nabla_{\bar k}u_j=g_{j\bar k}$, differentiating $|\partial u|_g^2$ gives
\begin{equation}\label{eq:bochner-gradient}
 \partial_i\bigl(|\partial u|_g^2\bigr)=u_i+g^{j\bar k}\mathcal H_{ij}u_{\bar k}.
\end{equation}
Fix a point and choose holomorphic normal coordinates for $g$ there.
Commuting covariant derivatives and using $\Delta_gu=n$ and $\Ric(g)=g$, we obtain
\[
 \sum_i\nabla_{\bar i}\mathcal H_{ij}=u_j.
\]
Differentiating \eqref{eq:bochner-gradient} once more therefore gives
\[
 \Delta_g\bigl(|\partial u|_g^2\bigr)=n+\sum_j u_ju_{\bar j}
             +\sum_{i,j}|\mathcal H_{ij}|^2
          =n+|\partial u|_g^2+|\mathcal H|_g^2.
\]
\end{proof}

The following differential inequality gives the nonnegative weighted divergence
used in the proof of Proposition~\ref{prop:normal}.
\begin{lemma}\label{lem:F-inequality}
Let $u$ be a smooth strictly plurisubharmonic solution of \eqref{eq:main}, and set
\begin{equation}\label{eq:F-definition}
 F=e^{-\frac{u}{n+1}}
       \left(n+1+|\partial u|_g^2\right).
\end{equation}
Then
\begin{equation}\label{eq:F-gradient}
 |\partial F|_g^2
 \leq\bigl[F-(n+1)\bigr]\Delta_gF.
\end{equation}
\end{lemma}
\begin{proof}
By the Bochner formula in Lemma~\ref{lem:bochner}, we obtain
\begin{eqnarray*}
 e^{\frac{u}{n+1}}\Delta_gF &=& \Delta_g\bigl(|\partial u|_g^2\bigr) -\frac{2}{n+1}\operatorname{Re} \left\langle \partial u,\partial\bigl(|\partial u|_g^2\bigr) \right\rangle_g\\
 &&+\left[\frac{|\partial u|_g^2}{(n+1)^2} -\frac{n}{n+1} \right] \left(n+1+|\partial u|_g^2\right)\\
  &=& |\mathcal H|_g^2 -\frac{2}{n+1}\operatorname{Re} \left\langle \mathcal H, \partial u\otimes\partial u \right\rangle_g +\frac{|\partial u|_g^4}{(n+1)^2}.
\end{eqnarray*}
Thus
\begin{equation}\label{eq:F-square}
 \Delta_gF
 =e^{-\frac{u}{n+1}}
   \left|
     \mathcal H
     -\frac{\partial u\otimes\partial u}{n+1}
   \right|_g^2
 \geq0.
\end{equation}
Similarly, \eqref{eq:bochner-gradient} gives
\[
 F_j
 =e^{-\frac{u}{n+1}}
   \sum_{k,\ell}g^{k\bar\ell}
   \left(\mathcal H_{jk}-\frac{u_ju_k}{n+1}\right)
   u_{\bar\ell}.
\]
By Cauchy--Schwarz and \eqref{eq:F-square},
\[
 |\partial F|_g^2
 \leq e^{-\frac{u}{n+1}}|\partial u|_g^2\Delta_gF.
\]
Since $u\leq0$,
\[
 F-(n+1)
 =e^{-\frac{u}{n+1}}|\partial u|_g^2
  +(n+1)\left(e^{-\frac{u}{n+1}}-1\right)
 \geq e^{-\frac{u}{n+1}}|\partial u|_g^2.
\]
Combining the last two inequalities gives \eqref{eq:F-gradient}.
\end{proof}

\subsection{Boundary geometry}\label{subsec:moments}

Denote by $\nu$ the Euclidean outward unit normal to $\partial B$
and by $d\sigma$ the Euclidean surface measure on $\partial B$.
The contact form and the complex tangent bundle on $\partial B$ are
\[
 \eta=\operatorname{Im}\sum_j\bar z_j\,dz_j,
 \qquad H=\ker\eta\subset T(\partial B).
\]
The circle action $z\mapsto e^{it}z$ generates the unit tangent field
$\xi(z)=\left.\frac{d}{dt}\right|_{t=0}e^{it}z=iz$ on $\partial B$.
In real coordinates $z_j=x_j+iy_j$, this field is
\[
 \xi=\sum_{j=1}^n\left(x_j\frac{\partial}{\partial y_j}
                      -y_j\frac{\partial}{\partial x_j}\right).
\]
It is the Reeb field of $\eta$ since $\eta(\xi)=|z|^2=1$ and
$d\eta(\xi,X)=-d(|z|^2)(X)=0$ for every $X\in T_z(\partial B)$.
If $\nabla_{\partial B}$ and $\Delta_{\partial B}$ are the Euclidean sphere gradient and Laplacian, set
\[
 \nabla_H f=\nabla_{\partial B} f-(\xi f)\xi,
 \qquad \Delta_H f=\Delta_{\partial B} f-\xi^2f.
\]
We use $H$ because it consists of the tangent directions whose multiples by $i$ are still tangent to $\partial B$. The Reeb direction lies outside $H$, since $i\xi=-\nu$ is normal to the sphere.

Since $\xi$ is divergence free on $\partial B$,
\begin{equation}\label{eq:H-ibp}
 \int_{\partial B} f\Delta_H g\dsigma
 =-\int_{\partial B}\langle\nabla_H f,\nabla_H g\rangle\dsigma.
\end{equation}
\subsection{Mass identity and Pohozaev identity}\label{subsec:identity}
We now express the Monge--Amp\`ere mass and Pohozaev identity in terms
of the boundary normal derivative.

Let $u\in C^\infty(\overline B)$ be strictly plurisubharmonic in $B$,
with $u=0$ on $\partial B$. Set $\rho=|z|^2-1$. Smooth division at the boundary gives
\[
 u=\rho h,\qquad h\in C^\infty\text{ near }\partial B,\qquad
 c:=h|_{\partial B}=\frac{1}{2}\partial_\nu u>0.
\]
The positivity follows from the Hopf lemma.
Here $c$ is a positive function on $\partial B$,
whose constancy for solutions of \eqref{eq:main} will be proved
in Section~\ref{sec:normal}.

The normal derivative determines the total Monge--Amp\`ere mass of $u$ on the unit ball.

\begin{lemma}\label{lem:boundary-mass}
Let $u$ be smooth and strictly plurisubharmonic on $\overline B$, with
$u=0$ on $\partial B$. Then
\begin{equation}\label{eq:boundary-mass}
 \int_B(\ddc u)^n=(2\pi)^n\fint_{\partial B}c^n\dsigma,
\end{equation}
where $\fint_{\partial B}$ denotes the surface average with respect to $d\sigma$.
For a solution of \eqref{eq:main}, the left-hand side equals $\Lambda$.
\end{lemma}
\begin{proof}
Let $\iota:\partial B\hookrightarrow\overline B$ denote the inclusion.
The contact form $\eta$ defined above satisfies
$\iota^*d^c\rho=\eta$ and $\iota^*(\ddc\rho)=d\eta$.
Since $u=\rho h$ and $\rho=0$ on $\partial B$, we have
\[
 \iota^*d^cu=c\eta,\qquad
 \iota^*(\ddc u)=d(c\eta)=dc\wedge\eta+c\,d\eta.
\]
Wedging with $c\eta$ eliminates every term containing $dc\wedge\eta$, so
\[
 \iota^*\!\left(d^cu\wedge(\ddc u)^{n-1}\right)
 =c^n\eta\wedge(d\eta)^{n-1}.
\]
Write the top-degree form $\eta\wedge(d\eta)^{n-1}$ as $a\,d\sigma$.
For every $U\in U(n)$, holomorphicity and $\rho\circ U=\rho$ give
$(U|_{\partial B})^*\eta=\iota^*d^c(\rho\circ U)=\eta$.
As $U$ also preserves $d\sigma$, it follows that $a(Uz)=a(z)$.
Thus $a$ is constant by transitivity of the $U(n)$-action on $\partial B$.

Since $\rho_{j\bar k}=\delta_{jk}$, we have $(\ddc\rho)^n=2^n n!\,dV$.
Stokes' theorem applied to $\rho$ gives
\[
 a|\partial B|=\int_{\partial B}\eta\wedge(d\eta)^{n-1}
 =\int_B(\ddc\rho)^n=2^n n!|B|=(2\pi)^n.
\]
Applying Stokes' theorem to $u$ now yields
\[
 \int_B(\ddc u)^n
 =a\int_{\partial B}c^n\dsigma
 =(2\pi)^n\fint_{\partial B}c^n\dsigma,
\]
as claimed.
\end{proof}

We use Chi Li's Pohozaev identity and sharp mass bound \cite{L16} for
$k=n$ on the unit ball. His parameter $a$ in the determinant equation
equals $\frac{\Lambda}{2^n n!}$. With the boundary mass formula
\eqref{eq:boundary-mass}, his result takes the following form.
\begin{lemma}\label{lem:pohozaev}
Every smooth solution of \eqref{eq:main} satisfies
\begin{equation}\label{eq:pohozaev}
 \fint_{\partial B}c^n\dsigma
   -\frac{1}{n+1}\fint_{\partial B}c^{n+1}\dsigma=\frac{\Lambda}{2^n n!Z}.
\end{equation}
In particular, $\Lambda<[2\pi(n+1)]^n$.
\end{lemma}

\subsection{Radial solutions}\label{subsec:radial}
We recall the classical classification of radial solutions \cite{D85,L16}.

\begin{proposition}\label{prop:radial-solutions}
The smooth radial plurisubharmonic solutions of \eqref{eq:main} on the unit ball are precisely
\begin{equation}\label{eq:radial-family}
 u(z)=(n+1)\log\frac{1+\mu|z|^2}{1+\mu},\qquad \mu=\frac{\Lambda^{\frac{1}{n}}}{2\pi(n+1)-\Lambda^{\frac{1}{n}}}>0.
\end{equation}
The solution with parameter $\mu$ has mass
\begin{equation}\label{eq:radial-mass}
 \Lambda=\left(\frac{2\pi(n+1)\mu}{1+\mu}\right)^n.
\end{equation}
In particular, for each $0<\Lambda<[2\pi(n+1)]^n$ there is a unique radial solution.
\end{proposition}
\begin{proof}
For a radial function $u(z)=\varphi(s)$, where $s=|z|^2$, differentiation gives
\[
 u_{j\bar k}=\varphi'(s)\delta_{jk}+\varphi''(s)\bar z_jz_k.
\]
The tangential and radial eigenvalues are respectively $\varphi'$ with multiplicity $n-1$ and $\varphi'+s\varphi''$ with multiplicity one. Thus
\[
 \det g=(\varphi')^{n-1}(\varphi'+s\varphi'').
\]

For a radial solution of \eqref{eq:main}, the equation becomes
\[
 (\varphi')^{n-1}(\varphi'+s\varphi'')=\frac{\Lambda}{2^n n!Z}e^{-\varphi},
 \qquad \varphi'>0,
\]
where for $n=1$ the positivity of $\varphi'$ follows by integrating
$(s\varphi')'>0$ from the origin. Multiplication by $ns^{n-1}$ and integration give
\[s^n \varphi'(s)^n=\frac{\Lambda}{2^n(n-1)!Z}\int_0^s t^{n-1}e^{-\varphi(t)}\,dt,\]
that is,
\begin{equation}\label{eq:radial-integral}
 \varphi'(s)^n=\frac{\Lambda}{2^n(n-1)!Z}\int_0^1 t^{n-1}e^{-\varphi(st)}\,dt.
\end{equation}

Set $\mu=\frac{\varphi'(0)}{n+1}$ and
$\varphi_*(s)=\varphi(0)+(n+1)\log(1+\mu s)$. Since
\[\varphi'(0)^n=\frac{\Lambda}{2^n n!Z}e^{-\varphi(0)},\] direct substitution shows that $\varphi_*$ solves the same integral equation. Both derivatives have positive lower bounds on $[0,1]$. The mean value theorem applied to the exponential and to the $n$th root in \eqref{eq:radial-integral} yields
\[
 |\varphi'(s)-\varphi_*'(s)|\leq C_0\sup_{0\leq t\leq s}|\varphi(t)-\varphi_*(t)|.
\]
If $D(s)=\displaystyle\sup_{0\leq t\leq s}|\varphi(t)-\varphi_*(t)|$, integration gives
\[\displaystyle D(s)\leq C_0\int_0^sD(t)\,dt.\] 
Gronwall's inequality implies $\varphi=\varphi_*$. Since $\varphi(1)=0$, this gives \eqref{eq:radial-family}.
For this function, $c=\frac{(n+1)\mu}{1+\mu}$, so \eqref{eq:boundary-mass} gives \eqref{eq:radial-mass}.
The mass is strictly increasing from $0$ to $[2\pi(n+1)]^n$ as $\mu$ ranges from $0$ to $\infty$, proving uniqueness within the radial class.

Direct substitution gives
\[
 \det g=\frac{((n+1)\mu)^n}{(1+\mu|z|^2)^{n+1}},
 \qquad Z=|B|(1+\mu).
\]
Together with \eqref{eq:radial-mass}, these identities verify
\eqref{eq:main} and give existence throughout \eqref{eq:range}.
\end{proof}

\begin{remark}
In the original normalization of \cite{BB22}, the equation is
\begin{equation}\label{eq:original}
 \left(\frac{i}{2\pi}\partial\pbar\varphi\right)^n
 =\frac{e^{-\gamma\varphi}}{\int_B e^{-\gamma\varphi}\dV}\,dV,
 \qquad \varphi|_{\partial B}=0.
\end{equation}
The substitution $u=\gamma\varphi$ gives \eqref{eq:main} with $\Lambda=(2\pi\gamma)^n$. Accordingly, the radial solution in this normalization is
\[
 \varphi_\gamma(z)=\frac{n+1}{\gamma}
 \log\frac{1+\mu_\gamma|z|^2}{1+\mu_\gamma},
 \qquad \mu_\gamma=\frac{\gamma}{n+1-\gamma},\qquad 0<\gamma<n+1.
\]
\end{remark}

\section{Smoothness of solutions}\label{sec:smoothness}
The subcritical maximizers of $\mathcal G_\Lambda$ belong to $\mathcal E^1(B)$. To apply the classification of smooth solutions to these maximizers, we first prove smoothness up to the boundary for bounded weak solutions and then show that every finite-energy solution is bounded.

\subsection{Smoothness by a Monge--Amp\`ere flow}
The smoothing follows the principle of Sz\'ekelyhidi and Tosatti \cite{ST11}. For the boundary
problem, we use Do's smoothing theorem and the uniqueness theorem of
Guedj, Lu and Zeriahi. We state them in the form needed below.

Let $T>0$, let $a\in C^\infty([0,T))$, and let
$u_0\in\PSH(B)\cap C(\overline B)$ with $u_0=0$ on $\partial B$.
Consider the problem
\[
 \begin{cases}
 \partial_t V=\log\det(V_{j\bar k})+a(t)
       &\text{in }(0,T)\times B,\\
 V(t,\zeta)=0&\zeta\in\partial B,\\
 V(0,z)=u_0(z)&z\in\overline B.
 \end{cases}
\]
Here a pluripotential solution has $V(t,\cdot)\in\PSH(B)$ for each $t$,
is locally uniformly Lipschitz in time, and satisfies
\[
 (\ddc V(t,\cdot))^n
 =2^n n!\,e^{\partial_t V-a(t)}\,dV
\]
in the Bedford--Taylor sense for almost every $t\in(0,T)$.

\begin{theorem}[{\cite[Theorem~0.1]{D17}}]
If $u_0$ extends to a bounded continuous plurisubharmonic function on
a neighborhood of $\overline B$, then the above problem has a solution
\[
 V\in C^\infty\bigl((0,T)\times\overline B\bigr)
       \cap C\bigl([0,T)\times\overline B\bigr).
\]
Moreover, $V(t,\cdot)$ is strictly plurisubharmonic in $B$ for every
$0<t<T$.
\end{theorem}

\begin{theorem}[{\cite[Corollary~6.6]{GLZ21}}]
Let $V,W\in C([0,T)\times\overline B)$ be two pluripotential solutions
of the above problem. If both are locally uniformly semiconcave in
$t\in(0,T)$, then $V=W$.
\end{theorem}

Here local uniform semiconcavity means that on every interval
$[\delta,S]\subset(0,T)$, the second time derivative is bounded above
in the sense of distributions, with a bound independent of $z\in B$.
The boundary data are identically zero, so their time derivatives vanish
and condition \textup{(6.4)} in \cite{GLZ21} is satisfied.
Since $a$ is smooth, the assumptions on the coefficients also hold on
each interval $[0,S]$ with $S<T$.

We now apply these results to the mean field equation.
\begin{proposition}\label{prop:smooth}
Let $n\geq1$. Every bounded weak solution of \eqref{eq:main} belongs to $C^\infty(\overline B)$. Its complex Hessian is uniformly positive definite on $\overline B$.
\end{proposition}
\begin{proof}
We compare the explicit pluripotential solution $W(t,z)=(1-t)u(z)$
with a solution of the same parabolic problem obtained by smoothing.
Uniqueness then gives smoothness of $u$.

Fix the given solution $u$. Its equation is
\[
 \frac{(\ddc u)^n}{2^n n!}=\frac{\Lambda}{2^n n!Z}e^{-u}\,dV.
\]
Since $u$ is bounded, this density is bounded above and below by positive constants. Ko\l odziej's Dirichlet theorem \cite{K98}, also recalled in \cite{GKZ08}, gives a continuous plurisubharmonic function on $\overline B$ with this prescribed density and zero boundary value. The Bedford--Taylor comparison principle \cite{BT82} identifies it with $u$. Thus $u\in C(\overline B)$.

To apply Do's smoothing theorem,
we first extend $u$ to a bounded continuous plurisubharmonic
function on a neighborhood of $\overline B$. 
Choose $b_0>0$ such that 
\[b_0^n>\frac{\Lambda}{2^n n!Z}\|e^{-u}\|_{L^\infty(B)}.\]
Since $\rho=|z|^2-1$ has identity complex Hessian, the determinant
density of $b_0\rho$ is $b_0^n$. The comparison principle and the
maximum principle therefore give
\begin{equation}\label{eq:barrier}
 b_0\rho\le u\le0\quad\hbox{in }B.
\end{equation}
Choose a sequence $\varepsilon_j>0$ decreasing to zero. For each $j$, define on the ball $B_{1+\delta}$, with any fixed $\delta>0$,
\[
 \widetilde u_j(z)=
 \begin{cases}
 \max\{u(z),b_0\rho(z)+\varepsilon_j\},& |z|<1,\\
 b_0\rho(z)+\varepsilon_j,&1\leq |z|<1+\delta.
 \end{cases}
\]
For each $j$, continuity of $u-b_0\rho$ and its vanishing on $\partial B$ imply that the two definitions agree with $b_0\rho+\varepsilon_j$ on a collar of $\partial B$. Hence $\widetilde u_j$ is plurisubharmonic. By \eqref{eq:barrier}, these functions decrease to
\[
 \widetilde u(z)=
 \begin{cases}
 u(z),&|z|\leq1,\\
 b_0\rho(z),&1<|z|<1+\delta.
 \end{cases}
\]
The limit is a bounded continuous plurisubharmonic extension of $u$.

For $0\leq t<\frac{1}{2}$, set $W(t,z)=(1-t)u(z)$.
Since $\partial_tW=-u$, homogeneity of the Bedford--Taylor operator gives
\[
 \frac{\bigl(\ddc W(t,\cdot)\bigr)^n}{2^n n!}
 =\frac{\Lambda(1-t)^n}{2^n n!Z}e^{\partial_tW}\,dV.
\]
Thus $W$ solves the parabolic problem above with
\[
 T=\frac{1}{2},\qquad u_0=u,\qquad
 a(t)=\log\frac{2^n n!Z}{\Lambda}-n\log(1-t).
\]
The function $a$ is smooth. The initial value has the bounded continuous
plurisubharmonic extension just constructed, and the zero boundary data
agree with $u$ on $\partial B$. Do's theorem therefore gives a solution
\[
 V\in C^\infty\bigl((0,\frac{1}{2})\times\overline B\bigr)
       \cap C\bigl([0,\frac{1}{2})\times \overline B\bigr),
\]
such that $V(t,\cdot)$ is strictly plurisubharmonic in $B$
for every $0<t<\frac{1}{2}$.

Both $V$ and $W$ are locally uniformly Lipschitz and semiconcave in time.
Indeed, $\partial_tW=-u$ is bounded and $\partial_t^2W=0$.
For every $0<\delta<S<\frac{1}{2}$, the derivatives
$\partial_tV$ and $\partial_t^2V$ are bounded on
$[\delta,S]\times\overline B$ by the smoothness of $V$.
The two solutions have the same initial value $u$ and zero boundary data.
Applying the uniqueness theorem above on $(0,S)\times B$ for each
$0<S<\frac{1}{2}$, we obtain
\[
 V(t,z)=W(t,z)=(1-t)u(z),\qquad 0<t<\frac{1}{2}.
\]
At $t=\frac{1}{4}$ this gives
\[
 u(z)=\frac{4}{3} V(\frac{1}{4},z)\in C^\infty(\overline B).
\]
The classical equation extends to $\partial B$ with positive density.
Thus $g\ge0$ and $\det g>0$ on the compact set $\overline B$, giving
uniform positive definiteness.
\end{proof}

\subsection{Finite-energy solutions}
It remains to reduce finite-energy solutions to bounded ones.
\begin{proposition}\label{prop:finite-energy}
Every solution of \eqref{eq:main} in $\mathcal E^1(B)$ is bounded and smooth on $\overline B$.
\end{proposition}
\begin{proof}
The Moser--Trudinger inequality of \cite[Theorem~1.5 and Section~4]{BB22}, extended by decreasing approximation, gives constants $A_n,B_n$ such that
\[
 \log\int_B e^{-w}\dV\leq A_n\mathcal E(w)+B_n,
 \qquad w\in\mathcal E^1(B).
\]
Apply this to $w=pu$. Since $\mathcal E(pu)=p^{n+1}\mathcal E(u)$, we obtain
\[
 \log\int_Be^{-pu}\dV\leq A_n p^{n+1}\mathcal E(u)+B_n<\infty
 \qquad(p>1).
\]
Thus the fixed density $\frac{\Lambda e^{-u}}{\int_Be^{-u}\dV}$ belongs to $L^p(B)$ for every finite $p>1$. Ko\l odziej's theorem gives a continuous bounded zero-boundary solution $v$ with the same Monge--Amp\`ere measure, and $v\in\mathcal E_0(B)\subset\mathcal E^1(B)$.

For $n\geq2$, the mixed-energy inequality \cite[Theorem~2.1]{ACC12} gives, for every $\psi\in\mathcal E^1(B)$,
\[
 \int_B(-\psi)(\ddc u)^n
 \leq
 (n+1)\mathcal E(\psi)^{\frac{1}{n+1}}\mathcal E(u)^{\frac{n}{n+1}}.
\]
Thus the Monge--Amp\`ere measure $(\ddc u)^n$, which also equals
$(\ddc v)^n$, satisfies the condition for uniqueness in
\cite[Proposition~3.1]{ACC12}. Since both functions belong to $\mathcal E^1(B)$, we obtain $u=v$. For $n=1$, zero-boundary finite-energy potentials belong to
$W^{1,2}_0(B)$, with
$\mathcal E(w)=\frac{1}{4}\displaystyle\int_B|\nabla w|^2\dV$.
Thus $u-v\in W^{1,2}_0(B)$, and equality of their Laplacians gives
\[\displaystyle\int_B|\nabla(u-v)|^2\dV=0,\] so again $u=v$. Hence $u\in\mathcal E_0(B)$, and the previous proposition gives smoothness on $\overline B$.
\end{proof}

\section{Radial symmetry and classification on the ball}\label{sec:normal}

We first prove constancy of the boundary normal derivative, then use it
to establish radial symmetry and complete the proof of Theorem~\ref{thm:MT}.

\subsection{Constancy of the boundary normal derivative}
We use the notation of Sections~\ref{subsec:notation},~\ref{subsec:moments} and~\ref{subsec:identity}.

We use the following form of the divergence theorem in complex coordinates,
obtained by taking real parts in \cite[Lemma~2.2]{L16} on the unit ball.

\begin{lemma}\label{lem:complex-divergence}
Let $X=\displaystyle\sum_{j=1}^n X^j\frac{\partial}{\partial z_j}$ be a $C^1$ vector field
of type $(1,0)$ on $\overline B$. Then
\[
 \int_B\operatorname{Re}
       \left(\sum_j\partial_jX^j\right)\dV
 =
 \frac{1}{2}\int_{\partial B}
       \operatorname{Re}
       \left(\sum_jX^j\bar z_j\right)\dsigma.
\]
\end{lemma}
We now prove constancy of the boundary normal derivative.
\begin{proposition}\label{prop:normal}
Every smooth solution of \eqref{eq:main} has constant boundary normal derivative. More precisely,
\[
 \partial_\nu u\big|_{\partial B}=\frac{\Lambda^{\frac{1}{n}}}{\pi}.
\]
\end{proposition}

\begin{proof}
Assume first that $n\geq2$. To prove constancy of $\partial_\nu u=2c$,
it suffices to show that $\nabla_Hc=0$ on $\partial B$.
Indeed, under this condition, $dc$ vanishes on $H=\ker\eta$, so
$dc=(\xi c)\eta$.
Taking the exterior derivative gives
\[
 0=d(\xi c)\wedge\eta+(\xi c)d\eta.
\]
At every point choose a nonzero $X\in H$. Since $iX\in H$ and
$d\eta(X,iX)=2|X|^2$, evaluation on $(X,iX)$ gives
$2(\xi c)|X|^2=0$.
Thus $\xi c=0$ and $dc=0$. The sphere is connected, so $c$ is constant.

Since $c>0$, we are therefore reduced to proving
\begin{equation}\label{eq:horizontal-energy-zero}
 \int_{\partial B}c^{n-2}|\nabla_Hc|^2\dsigma=0.
\end{equation}
We prove \eqref{eq:horizontal-energy-zero} by integrating a nonnegative
weighted divergence and using the Pohozaev identity to cancel the terms
without derivatives of $c$.

\smallskip
Let $F$ be defined by \eqref{eq:F-definition}. 
Since $u<0$ in $B$ and $\partial_\nu u=2c>0$ on $\partial B$, we have $F-(n+1)>0$ on $\overline B$. Using Lemma~\ref{lem:F-inequality} and the cofactor identity \eqref{eq:cofactor}, we obtain
\begin{eqnarray}
 \lefteqn{\sum_{j,k}\partial_j \left[ \bigl(F-(n+1)\bigr)^{-\frac{n}{n+1}} \det g\,g^{j\bar k}F_{\bar k} \right]}\nonumber\\
  &=& \det g\cdot\bigl[F-(n+1)\bigr]^{-\frac{n}{n+1}} \left( \Delta_gF -\frac{n}{n+1} \frac{|\partial F|_g^2}{F-(n+1)} \right)\nonumber\\
  &\geq& \frac{\det g}{n+1} \bigl[F-(n+1)\bigr]^{-\frac{n}{n+1}} \Delta_gF \geq0. \label{eq:weighted-divergence}
\end{eqnarray}
Integrating \eqref{eq:weighted-divergence} over $B$ and applying
Lemma~\ref{lem:complex-divergence}, we obtain
\begin{equation}\label{bound-int}
 0\leq\frac{1}{2}\int_{\partial B}
 \bigl[F-(n+1)\bigr]^{-\frac{n}{n+1}}
 \operatorname{Re}\left(\sum_{j,k}\det g\,g^{j\bar k}F_{\bar k}\bar z_j\right)\dsigma.
\end{equation}
To compute this boundary integral, we next express the complex Hessian of $u$, the boundary value of $F$ and its normal derivative in terms of $c$.

\begin{itemize}
\item \textit{The complex Hessian on the boundary.}
Fix $P\in\partial B$ and make a unitary change of coordinates
so that $P=(0,\ldots,0,1)$.
Write $z=(w,z_n)$, where $w\in\C^{n-1}$, and use the boundary chart
\[
 (w,\theta)\longmapsto
 \bigl(w,e^{i\theta}\sqrt{1-|w|^2}\bigr).
\]
The derivatives $c_\alpha$ and $c_{\alpha\bar\beta}$ below
are derivatives of the pullback of $c$ to this chart at $(0,0)$.
At $P$,
\begin{equation}\label{eq:local-H}
 \sum_{\alpha<n}|c_\alpha|^2
 =\frac{1}{4}|\nabla_Hc|^2,
 \qquad
 \sum_{\alpha<n}c_{\alpha\bar\alpha}
 =\frac{1}{4}\Delta_Hc.
\end{equation}
In the ambient coordinates, write $z_n=r+iy$.
Thus $\partial_r=\partial_\nu$ at $P$.

Recall that $u=\rho h$, where $\rho=|z|^2-1$ and
$h|_{\partial B}=c$.
Differentiating this identity gives
\begin{eqnarray*}
 u_j&=&\rho_jh+\rho h_j,\\
 u_{j\bar k}&=&\delta_{jk}h+\rho_jh_{\bar k}
                  +\rho_{\bar k}h_j+\rho h_{j\bar k}.
\end{eqnarray*}
At $P$, we have $\rho=0$, $\rho_\alpha=0$ for $\alpha<n$,
and $\rho_n=1$. The first derivatives at $P$ are
\begin{equation}\label{eq:boundary-first}
 u_\alpha=0\quad(\alpha<n),\qquad u_n=c.
\end{equation}
The complex Hessian at $P$ is
\begin{equation}\label{eq:boundary-block}
 g(P)=
 \begin{pmatrix}
  cI_{n-1}&\ \ (c_\alpha)_{\alpha<n}\\[6pt]
  (\overline{c_\beta})_{\beta<n}&\ \ c+h_r
 \end{pmatrix}.
\end{equation}
In particular, the tangential block is $cI_{n-1}$,
because $\rho_{\alpha\bar\beta}=\delta_{\alpha\beta}$.
\vskip 10pt
\item \textit{Boundary value of $F$.}
Expanding the determinant in \eqref{eq:boundary-block} and using
\eqref{eq:local-H}, we obtain
\begin{equation}\label{eq:boundary-determinant}
 \det g=c^{n-1}(c+h_r)-c^{n-2}\sum_{\alpha<n}|c_\alpha|^2
       =c^{n-1}(c+h_r)-\frac{1}{4}c^{n-2}|\nabla_Hc|^2.
\end{equation}
The cofactor formula gives at $P$
\begin{equation}\label{eq:boundary-cofactors}
 \begin{aligned}
  \det g\cdot g^{n\bar n}&=c^{n-1},\\
  \det g \cdot g^{n\bar\alpha}&=-c^{n-2}c_\alpha,\qquad
  \det g\cdot g^{\alpha\bar n}=-c^{n-2}\overline{c_\alpha},
  \quad \alpha<n.
 \end{aligned}
\end{equation}
By \eqref{eq:boundary-first} and \eqref{eq:boundary-cofactors}, we have
\[
 |\partial u|_g^2
 =\sum_{j,k}g^{j\bar k}u_ju_{\bar k}
 =g^{n\bar n}u_nu_{\bar n}
 =c^2g^{n\bar n}
 =\frac{c^{n+1}}{\det g}.
\]
By \eqref{eq:det-equation} and $u=0$ on $\partial B$,
we have $\det g=\frac{\Lambda}{2^n n!Z}$ on $\partial B$.
Substituting into \eqref{eq:F-definition}, we obtain on $\partial B$
\begin{equation}\label{eq:F-boundary}
 F=n+1+\frac{2^n n!Z}{\Lambda}c^{n+1},
\end{equation}
that is,
\begin{equation}\label{eq:J-boundary}
 \bigl[F-(n+1)\bigr]^{-\frac{n}{n+1}}
 =\left(\frac{\Lambda}{2^n n!Z}\right)^{\frac{n}{n+1}}c^{-n}.
\end{equation}
\vskip 10pt
\item \textit{The normal derivative of $F$.}
Normal differentiation must use the interior definition of $F$.
All quantities below are evaluated at $P$.
Differentiating $u=\rho h$ in the $r$ direction gives
\[
 (u_\alpha)_r=2c_\alpha,\qquad
 \operatorname{Re}(u_n)_r=c+2h_r,\qquad
 u_r=2c.
\]
In the boundary chart, $r_{\alpha\bar\beta}=-\frac{1}{2}\delta_{\alpha\beta}$
at $P$. Hence
\[
 c_{\alpha\bar\beta}=h_{\alpha\bar\beta}
       -\frac{1}{2}\delta_{\alpha\beta}h_r,
 \qquad (u_{\alpha\bar\beta})_r
       =\delta_{\alpha\beta}h_r+2h_{\alpha\bar\beta}.
\]
Taking the trace and using \eqref{eq:local-H} gives
\[
 \sum_{\alpha<n}(u_{\alpha\bar\alpha})_r
 =2(n-1)h_r+\frac{1}{2}\Delta_Hc.
\]

In these fixed coordinates, the cofactor identity
\[
 g^{n\bar n}
 =\frac{\det(u_{\alpha\bar\beta})_{\alpha,\beta<n}}{\det g}
\]
holds in a neighborhood of $P$ and can therefore be differentiated
in the normal direction. Differentiating \eqref{eq:det-equation} at $P$
gives $\partial_r\log\det g=-u_r=-2c$.
Since $(u_{\alpha\bar\beta})=cI_{n-1}$ at $P$, we obtain
\begin{equation}\label{eq:inverse-normal}
  \frac{\partial_rg^{n\bar n}}{g^{n\bar n}}
  =\frac{1}{c}\sum_{\alpha<n}(u_{\alpha\bar\alpha})_r
    -\partial_r\log\det g
  =2c+\frac{2(n-1)h_r}{c}+\frac{\Delta_Hc}{2c}.
\end{equation}

Next, differentiate the interior expression
$|\partial u|_g^2=g^{j\bar k}u_ju_{\bar k}$.
Using \eqref{eq:boundary-first} and \eqref{eq:boundary-cofactors}, we have
\begin{eqnarray*}
 \partial_r\bigl(|\partial u|_g^2\bigr) &=& 2\operatorname{Re} \sum_{j,k}g^{j\bar k}(u_j)_r u_{\bar k} +\sum_{j,k}(\partial_rg^{j\bar k})u_ju_{\bar k}\nonumber\\
  &=& 2cg^{n\bar n}\operatorname{Re} \left( (u_n)_r-\frac{1}{c} \sum_{\alpha<n}\overline{c_\alpha}(u_\alpha)_r \right) +c^2\partial_rg^{n\bar n}\nonumber\\
  &=& 2cg^{n\bar n} \left(c+2h_r-\frac{|\nabla_Hc|^2}{2c}\right) +c^2\partial_rg^{n\bar n}\nonumber\\
  &=& 4c-2c^2g^{n\bar n}+c^2\partial_rg^{n\bar n}.
\end{eqnarray*}
In the last equality we used \eqref{eq:boundary-cofactors} and
\[
 h_r=\frac{\det g}{c^{n-1}}-c+\frac{|\nabla_Hc|^2}{4c},
\]
which follows from \eqref{eq:boundary-determinant}.

Differentiating \eqref{eq:F-definition} at $P$, we have
\[
 F_r=\partial_r\bigl(|\partial u|_g^2\bigr)
     -\frac{2c}{n+1}\bigl(n+1+|\partial u|_g^2\bigr).
\]
Substituting \eqref{eq:inverse-normal} and the above expression for $h_r$,
and using \eqref{eq:boundary-cofactors}, we obtain
\begin{equation}\label{eq:F-normal}
  \det g\cdot F_r=2nc\,\det g-2nc^{n+1}
       +\frac{2n}{n+1}c^{n+2}+\frac{n-1}{2}c^{n-1}|\nabla_Hc|^2
       +\frac{1}{2}c^n\Delta_Hc.
\end{equation}

\end{itemize}

Now we return to \eqref{bound-int}.
Differentiating \eqref{eq:F-boundary} in the complex tangential directions
and using the boundary value of $\det g$, we obtain
\[
 F_\alpha
 =\frac{(n+1)c^n}{\det g}c_\alpha.
\]
Using  \eqref{eq:local-H}, \eqref{eq:boundary-cofactors} and
$F_{\bar n}=\frac{F_r+iF_y}{2}$, we find
\begin{eqnarray*}
 \frac{1}{2}\operatorname{Re}
       \sum_k\det g\,g^{n\bar k}F_{\bar k}
 &=& \frac{c^{n-1}}{2}\operatorname{Re}
 \left(
   F_{\bar n}
   -\frac{1}{c}\sum_{\alpha<n}c_\alpha F_{\bar\alpha}
 \right)
 \notag\\
 &=& \frac{c^{n-1}}{4}
 \left(
   F_r-\frac{2(n+1)c^{n-1}}{\det g}
             \sum_{\alpha<n}|c_\alpha|^2
 \right)
 \notag\\
 &=&
 \frac{c^{n-1}}{4\det g}
 \left(
   \det g \cdot F_r-\frac{n+1}{2}c^{n-1}|\nabla_Hc|^2
 \right).
\end{eqnarray*}
Substituting \eqref{eq:F-normal} and using
$\det g=\frac{\Lambda}{2^n n!Z}$ on $\partial B$ gives
\begin{equation}\label{eq:expanded-flux}
\begin{aligned}
 \frac{1}{2}\operatorname{Re}\sum_k\det g\,g^{n\bar k}F_{\bar k}
 &=\frac{2^n n!Z}{4\Lambda}\Bigg[
 2nc^n\left(\frac{\Lambda}{2^n n!Z}-c^n+\frac{c^{n+1}}{n+1}\right)\\
 &\qquad{}-c^{2n-2}|\nabla_Hc|^2+\frac{c^{2n-1}}{2}\Delta_Hc\Bigg].
\end{aligned}
\end{equation}
In unitary coordinates with $z=(0,\ldots,0,1)$ at the boundary point,
only the term $j=n$ contributes to the boundary integral in \eqref{bound-int}.
Substituting  \eqref{eq:J-boundary} and \eqref{eq:expanded-flux}
into this integral, we obtain
\[
 0\leq\frac{1}{4}\left(\frac{2^n n!Z}{\Lambda}\right)^{\frac{1}{n+1}}
 \int_{\partial B}\Bigg[
 2n\left(\frac{\Lambda}{2^n n!Z}-c^n+\frac{c^{n+1}}{n+1}\right)-c^{n-2}|\nabla_Hc|^2
 +\frac{1}{2}c^{n-1}\Delta_Hc\Bigg]\dsigma.
\]
Dividing by the positive prefactor and integrating term by term gives
\begin{eqnarray}
 0 &\leq& \frac{2n\Lambda|\partial B|}{2^n n!Z} -2n\int_{\partial B}c^n\dsigma +\frac{2n}{n+1} \int_{\partial B}c^{n+1}\dsigma\nonumber\\
 && {}-\int_{\partial B} c^{n-2}|\nabla_Hc|^2\dsigma +\frac{1}{2}\int_{\partial B} c^{n-1}\Delta_Hc\dsigma. \label{eq:flux-integral}
\end{eqnarray}
The Pohozaev identity \eqref{eq:pohozaev} gives
\[
 \frac{2n\Lambda|\partial B|}{2^n n!Z}
 -2n\int_{\partial B}c^n\dsigma
 +\frac{2n}{n+1}\int_{\partial B}c^{n+1}\dsigma=0.
\]
Moreover, \eqref{eq:H-ibp} gives
\[
 \int_{\partial B}c^{n-1}\Delta_Hc\dsigma
 =-(n-1)\int_{\partial B}c^{n-2}|\nabla_Hc|^2\dsigma.
\]
Substituting these identities into \eqref{eq:flux-integral}, we obtain
\[
 0\leq-\frac{n+1}{2}\int_{\partial B}c^{n-2}|\nabla_Hc|^2\dsigma\leq0.
\]
This proves \eqref{eq:horizontal-energy-zero}. Since $c>0$, we have
$\nabla_Hc=0$ on $\partial B$, and the argument at the beginning of the proof
shows that $c$ is constant. The mass identity \eqref{eq:boundary-mass} gives
\[
 c=\frac{\Lambda^{\frac{1}{n}}}{2\pi},\qquad
 \partial_\nu u=2c=\frac{\Lambda^{\frac{1}{n}}}{\pi}.
\]

For $n=1$, apply \cite{GNN79} to $w=-u>0$, which satisfies
\[
 -\Delta_{\mathrm{Eucl}}w=\frac{2\Lambda}{Z}e^w\quad\text{in }B\subset\R^2,
 \qquad w|_{\partial B}=0.
\]
The moving-plane theorem makes $u$ radial, so its normal derivative
is constant. The mass identity \eqref{eq:boundary-mass} then gives
$\partial_\nu u=\frac{\Lambda}{\pi}$.
\end{proof}

\subsection{Classification of solutions}
We now complete the proof of Theorem~\ref{thm:MT}.

\begin{proof}[Proof of Theorem~\ref{thm:MT}]
Let $u$ be a smooth solution of \eqref{eq:main}.
By Proposition~\ref{prop:normal}, $\partial_\nu u$ is constant on $\partial B$.
We show that $u$ is radial.
Let $K$ be a skew-Hermitian matrix, put $U_s=e^{sK}$, and define
\[
 u_s(z)=u(U_sz),\qquad
 v(z)=\left.\frac{d}{ds}\right|_{s=0}u_s(z).
\]
Unitary changes of variables preserve $B$ and Euclidean volume, so $\displaystyle\int_B e^{-u_s}\dV=Z$. Thus every $u_s$ satisfies
\[\det((u_s)_{j\bar k})=\frac{\Lambda}{2^n n!Z}e^{-u_s}.\]
Differentiating at $s=0$ gives
\[
 \Delta_gv+v=0\qquad\text{in }B.
\]
Since $u_s=0$ on $\partial B$, we have $v=0$ there. Moreover,
$U_s\nu(z)=\nu(U_sz)$, so
\[
 \partial_\nu u_s(z)=\partial_\nu u(U_sz)
\]
for $z\in\partial B$.
The right-hand side is independent of $s$, and therefore
$\partial_\nu v=0$ on $\partial B$.

The matrix $g$ is smooth and uniformly positive definite on
$\overline B$. In real coordinates, $\Delta_g+1$ is therefore a
uniformly elliptic operator with smooth coefficients. By boundary unique continuation for second-order elliptic equations
\cite[Remark~2]{Ar57}, the conditions
$v=\partial_\nu v=0$ on $\partial B$ imply $v=0$ in $B$.
Indeed, extend the coefficients smoothly across the boundary,
preserving ellipticity, and extend $v$ by zero. Its zero value and
zero normal derivative make this an $H^2_{\mathrm{loc}}$ solution
of the extended equation. It vanishes on an exterior open set,
so unique continuation gives the conclusion.

Since $K$ was arbitrary, every infinitesimal unitary rotation
annihilates $u$. The group $U(n)$ is connected and acts transitively
on each sphere centered at the origin, so $u$ is radial.

Proposition~\ref{prop:radial-solutions} now gives the formula \eqref{eq:explicit},
the necessary mass restriction \eqref{eq:range}, and uniqueness.
Conversely, for every $\Lambda$ satisfying \eqref{eq:range},
the radial family \eqref{eq:radial-family} provides a smooth solution.
Propositions~\ref{prop:smooth} and~\ref{prop:finite-energy} give the corresponding
assertions for bounded and finite-energy solutions.
\end{proof}

\begin{remark}
For a general positive smooth density $f$, it is natural to ask which
conditions on $f$ allow the equation $\det(u_{j\bar k})=f(u)$ to admit
an auxiliary function with a signed divergence and yield constancy of
the boundary normal derivative.
\end{remark}

\section{The sharp logarithmic inequality on the ball}\label{sec:sharp-log-ball}

We now use the classification in Theorem~\ref{thm:MT} to identify the subcritical maximizers of $\mathcal G_\Lambda$. Their explicit form gives both optimal constants in Theorem~\ref{thm:sharp-log}.

\begin{proof}[Proof of Theorem~\ref{thm:sharp-log}]
We follow the argument after Conjecture~7.6 in \cite{BB22}.
We first prove the inequality by passing to the limit in the subcritical range.
For $\mu>0$, let $u_{\mu}$ be the radial function in \eqref{eq:radial-family}, and set
\[
 \Lambda_{\mu}=\left(\frac{2\pi(n+1)\mu}{1+\mu}\right)^n.
\]
The variational existence theorem \cite[Theorem~7.2]{BB22}, with
$\gamma=\frac{(n+1)\mu}{1+\mu}$ and the change of normalization in
\eqref{eq:original}, gives a maximizer of $\mathcal G_{\Lambda_{\mu}}$
on $\mathcal E^1(B)$ which solves \eqref{eq:main} with mass $\Lambda_{\mu}$.
This existence result uses only the quasi-sharp inequality in that paper.
By Proposition~\ref{prop:finite-energy} and Theorem~\ref{thm:MT},
the maximizer must be $u_{\mu}$.
Dividing the maximality inequality by $\Lambda_{\mu}$ and subtracting
$\log|B|$, we obtain, for every $v\in\mathcal E^1(B)$,
\begin{equation}\label{eq:mt-subcritical}
 \log\left(\frac{1}{|B|}\int_Be^{-v}\dV\right)
   -\frac{\mathcal E(v)}{\Lambda_{\mu}}\leq
   \log\left(\frac{1}{|B|}\int_Be^{-u_{\mu}}\dV\right)
   -\frac{\mathcal E(u_{\mu})}{\Lambda_{\mu}}.
\end{equation}

Section~\ref{subsec:radial} gives
\[\displaystyle\int_Be^{-u_{\mu}}\dV=|B|(1+\mu).\]
Radial integration and integration by parts also give
\begin{eqnarray}
 \mathcal E(u_{\mu}) &=& [2\pi(n+1)]^n\int_0^\mu\frac{t^n}{(1+t)^{n+1}}\,dt\nonumber\\
  &=& [2\pi(n+1)]^n \left[\log(1+\mu)-\sum_{j=1}^n\frac{1}{j} \left(\frac{\mu}{1+\mu}\right)^j\right]. \label{eq:mt-radial-energy}
\end{eqnarray}
Since $\Lambda_\mu<[2\pi(n+1)]^n$ and $\mathcal E(u_\mu)\geq0$, \eqref{eq:mt-radial-energy} gives
\begin{eqnarray*}
 \log\left(\frac{1}{|B|}\int_B e^{-u_\mu}\dV\right)
 -\frac{\mathcal E(u_\mu)}{\Lambda_\mu}
 &\leq& \log(1+\mu)-\frac{\mathcal E(u_\mu)}{[2\pi(n+1)]^n}\\[5pt]
 &=& \sum_{j=1}^n\frac{1}{j}\left(\frac{\mu}{1+\mu}\right)^j
 \leq\sum_{j=1}^n\frac{1}{j}.
\end{eqnarray*}
For fixed $v\in\mathcal E^1(B)$, use this bound in \eqref{eq:mt-subcritical} and let $\mu\to\infty$. Since $\Lambda_\mu\to [2\pi(n+1)]^n$, we obtain \eqref{eq:sharp-mt}.

We next prove optimality of the additive constant for the stated energy coefficient.
By \eqref{eq:mt-radial-energy},
\[
 \log\left(\frac{1}{|B|}\int_Be^{-u_\mu}\dV\right)
 -\frac{\mathcal E(u_\mu)}{[2\pi(n+1)]^n}
 =\sum_{j=1}^n\frac{1}{j}\left(\frac{\mu}{1+\mu}\right)^j
 \longrightarrow\sum_{j=1}^n\frac{1}{j}
\]
as $\mu\to\infty$.
Thus the additive constant cannot be decreased.

Finally, \eqref{eq:mt-radial-energy} gives
\[\mathcal E(u_\mu)=[2\pi(n+1)]^n\log(1+\mu)+O(1),\]
while the logarithmic integral equals $\log(1+\mu)$.
Hence any smaller energy coefficient fails for this radial family,
even with an arbitrary finite additive constant.
\end{proof}

\begin{remark}
Theorem~\ref{thm:sharp-log} controls the functional on the whole class
$\mathcal E^1(B)$. Substituting \eqref{eq:sharp-mt} into the definition
of $\mathcal G_\Lambda$ gives, for every $u\in\mathcal E^1(B)$ and $\Lambda>0$,
\[
 \mathcal G_\Lambda(u)
 \leq-\left(1-\frac{\Lambda}{[2\pi(n+1)]^n}\right)
 \mathcal E(u)+\Lambda\left(\log|B|+\sum_{j=1}^n\frac{1}{j}\right).
\]
Thus $\mathcal G_\Lambda\to-\infty$ as $\mathcal E\to\infty$ in the
subcritical range, and remains bounded above at the critical mass.
For larger masses, the radial test functions used above satisfy
\[
 \mathcal G_\Lambda(u_{\mu})
 =\bigl(\Lambda-[2\pi(n+1)]^n\bigr)\log(1+\mu)+O(1)
 \longrightarrow+\infty
\]
as $\mu\to\infty$ by \eqref{eq:mt-radial-energy}.
\end{remark}

\section{Critical volume--capacity and  exponential Moser-Trudinger inequalities}\label{sec:exponential-ball}

We apply the logarithmic inequality to prove volume--capacity and  exponential 
Moser-Trudinger inequalities.

\subsection{Volume--capacity and sublevel estimates}

We first apply the logarithmic inequality to relative extremal functions to prove Theorem~\ref{thm:volume-capacity} and Corollary~\ref{cor:sublevel-volume}. 
\begin{proof}[Proof of Theorem~\ref{thm:volume-capacity}]
Let
\[
 h_K=\sup\{v\in\PSH(B):v\leq-1\text{ on }K,\ v\leq0\text{ on }B\}^*
\]
be the relative extremal function of $K$ in $B$. If $\Capa(K,B)=0$, then $|K|=0$ and the estimate is immediate. Otherwise, $h_K\in\mathcal E_0(B)$, $h_K=-1$ quasi-everywhere on $K$, and
\[
 \Capa(K,B)=\int_B(-h_K)(\ddc h_K)^n=(n+1)\mathcal E(h_K).
\]
Here quasi-everywhere means outside a pluripolar set, which has zero Lebesgue measure. Applying Theorem~\ref{thm:sharp-log} to $\lambda h_K$ and using
\[
 |K|e^\lambda\leq\int_Be^{-\lambda h_K}\dV,
\]
we obtain
\[
 \log\frac{|K|}{|B|}
 \leq\sum_{j=1}^n\frac{1}{j}
 +\frac{\lambda^{n+1}\Capa(K,B)}{(n+1)[2\pi(n+1)]^n}
 -\lambda.
\]
Taking $\lambda=\frac{2\pi(n+1)}{\Capa(K,B)^{\frac{1}{n}}}$ gives \eqref{eq:sharp-volume-capacity}.

For optimality, let $K=\overline{B_r}$, where $0<r<1$. Its relative extremal function is
\[
 h_K(z)=\max\left\{-1,\frac{\log|z|}{\log\frac{1}{r}}\right\},
\]
and hence
\[
 \Capa(\overline{B_r},B)=\frac{\pi^n}{\left(\log\frac{1}{r}\right)^n},
 \qquad
 \frac{|B_r|}{|B|}
 =\exp\left[-\frac{2\pi n}{\Capa(\overline{B_r},B)^{\frac{1}{n}}}\right].
\]
Letting $r\to0$ shows that the coefficient $2\pi n$ cannot be increased, even if the prefactor is replaced by an arbitrary finite constant.
\end{proof}

\begin{proof}[Proof of Corollary~\ref{cor:sublevel-volume}]
For $K\Subset\Omega\Subset B_R(a)$, monotonicity of relative capacity gives
\[
 \Capa(K,\Omega)\geq\Capa(K,B_R(a)).
\]
By translation and scaling, Theorem~\ref{thm:volume-capacity} therefore gives
\begin{equation}\label{eq:domain-volume-capacity}
 |K|\leq|B_R|\exp\left(\sum_{j=1}^n\frac{1}{j}
 -\frac{2\pi n}{\Capa(K,\Omega)^{\frac{1}{n}}}\right).
\end{equation}
The same estimate holds for Borel subsets of $\Omega$ by approximation. For $u\in\mathcal F(\Omega)$, the comparison principle gives
\[
 \Capa(\{u<-s\},\Omega)
 \leq\frac{1}{s^n}\int_\Omega(\ddc u)^n;
\]
see also \cite[Section~6.1]{ACKPZ09}. Substitution in \eqref{eq:domain-volume-capacity} proves \eqref{eq:sublevel-volume}.

The coefficient is optimal already for $\Omega=B$ and $u(z)=\log|z|$. In our normalization,
\[
 \int_B(\ddc u)^n=\pi^n,
 \qquad |\{u<-s\}|=|B|e^{-2ns}.
\]
Thus a larger coefficient fails as $s\to\infty$.
\end{proof}

\subsection{Exponential Moser-Trudinger inequality}
We then use the volume-capacity inequality to obtain the critical exponential inequality on the ball. 
It suffices to prove the following implication.
\begin{proposition}\label{prop:capacity-to-MT}
Let $\mu$ be a finite nonnegative Borel measure on $B$. Suppose that, for some $A,\beta>0$ and every compact set $K\Subset B$,
\begin{equation}\label{eq:volume-capacity}
 \mu(K)\leq A\exp\left[-\frac{\beta}{\Capa(K,B)^{\frac{1}{n}}}\right].
\end{equation} 
Then for every nonzero $u\in\mathcal H_0(B)$,
\[
 \int_B\exp\left[\beta\frac{(-u)^{\frac{n+1}{n}}}
 {[(n+1)\mathcal E(u)]^{\frac{1}{n}}}\right]d\mu
 \leq AC_n,
\]
where $C_n$ depends only on $n$.
\end{proposition}
Taking $d\mu=\frac{dV}{|B|}$, $A=\exp\left(\sum_{j=1}^n\frac{1}{j}\right)$ and $\beta=2\pi n$ gives Theorem~\ref{thm:sharp-MT} on the unit ball. Before proving the proposition, we recall the one-dimensional Moser lemma \cite{M71,BB14} and state the capacity comparison used below.

\begin{lemma}\label{lem:Moser}
For every integer $n\geq1$, there exists a constant $C_n>0$ with the following property. Let $w:(-\infty,0]\to\mathbb R$ be locally absolutely continuous and nondecreasing, with $w(0)=0$. If
\[
 \int_{-\infty}^0|w'(t)|^{n+1}\,dt\leq1,
\]
then
\[
 \int_{-\infty}^0\exp\left[(-w(t))^{\frac{n+1}{n}}\right]e^t\,dt\leq C_n.
\]
\end{lemma}

\begin{lemma}\label{lem:nested-capacity}
Let $F$ be compact and $E$ be a bounded hyperconvex open set with $F\Subset E\Subset B$. Then
\[
 \Capa(F,B)^{-\frac{1}{n}}
 \geq\Capa(\overline E,B)^{-\frac{1}{n}}
 +\Capa(F,E)^{-\frac{1}{n}}.
\]
\end{lemma}
\begin{proof}
The zero-capacity case is immediate, so assume that $\Capa(F,E)>0$.
Consider the relative extremal functions $h_{\overline E,B}$ and $h_{F,E}$, and extend $h_{F,E}$ by zero outside $E$.
Since $E$ is hyperconvex, $h_{F,E}$ has zero boundary limits, and its zero extension is quasi-continuous on $B$.
Fix $a>0$ and, for $b\geq0$, set
\[
 \psi_b=ah_{\overline E,B}+bh_{F,E},\qquad
 P_b=\sup\{v\in\PSH(B):v\leq\psi_b\text{ quasi-everywhere}\}^*.
\]
As $h_{\overline E,B}=-1$ in $E$, we have
\[
 (a+b)h_{\overline E,B}\leq P_b\leq0,
\]
so $P_b\in\mathcal E_0(B)$.
For $\varepsilon>0$, let
\[
 D_{b+\varepsilon}=\{P_{b+\varepsilon}=\psi_{b+\varepsilon}\}.
\]
Since $\psi_{b+\varepsilon}$ is quasi-continuous, \cite[Theorem~2.7]{ACLR25} gives
\[
 (\ddc P_{b+\varepsilon})^n
 =\mathbf1_{D_{b+\varepsilon}}(\ddc P_{b+\varepsilon})^n.
\]
On $D_{b+\varepsilon}$, outside a pluripolar set, we have
\[
 0\leq P_b-P_{b+\varepsilon}
 \leq\psi_b-\psi_{b+\varepsilon}
 =\varepsilon(-h_{F,E}).
\]
Moreover, $\psi_{b+\varepsilon}=-a+(b+\varepsilon)h_{F,E}$ in $E$.
Comparison on the contact set \cite[Theorem~2.2]{ACLR25} therefore gives
\[
 \mathbf1_{D_{b+\varepsilon}\cap E}(\ddc P_{b+\varepsilon})^n
 \leq(b+\varepsilon)^n\mathbf1_{D_{b+\varepsilon}\cap E}(\ddc h_{F,E})^n.
\]
By the convexity of $\mathcal E$ and the zero extension of $h_{F,E}$,
\[
\begin{aligned}
 \mathcal E(P_{b+\varepsilon})-\mathcal E(P_b)
 &\leq\int_B(P_b-P_{b+\varepsilon})(\ddc P_{b+\varepsilon})^n\\
 &\leq\varepsilon\int_{D_{b+\varepsilon}\cap E}(-h_{F,E})(\ddc P_{b+\varepsilon})^n\\
 &\leq\varepsilon(b+\varepsilon)^n\int_E(-h_{F,E})(\ddc h_{F,E})^n\\
 &=\varepsilon(b+\varepsilon)^n\Capa(F,E).
\end{aligned}
\]
Since $P_0=ah_{\overline E,B}$, summing this estimate over a partition of $[0,b]$ and letting the maximum length of the subintervals tend to zero gives
\[
 \mathcal E(P_b)
 \leq\frac{a^{n+1}}{n+1}\Capa(\overline E,B)
 +\frac{b^{n+1}}{n+1}\Capa(F,E).
\]
On the other hand, $P_b\leq-(a+b)$ quasi-everywhere on $F$, so
\[
 \frac{(a+b)^{n+1}}{n+1}\Capa(F,B)\leq\mathcal E(P_b).
\]
Take $a=\Capa(\overline E,B)^{-\frac{1}{n}}$ and $b=\Capa(F,E)^{-\frac{1}{n}}$. The upper bound for $\mathcal E(P_b)$ becomes $\frac{a+b}{n+1}$. Combining the two bounds gives
\[
 \Capa(F,B)^{-\frac{1}{n}}
 \geq a+b
 =\Capa(\overline E,B)^{-\frac{1}{n}}+\Capa(F,E)^{-\frac{1}{n}}.
\]
\end{proof}

\begin{proof}[Proof of Proposition~\ref{prop:capacity-to-MT}]
Let
\[
 E_s=\{u<-s\},\qquad F_s=\{u\leq-s\},\qquad
 m(s)=\int_{E_s}(\ddc u)^n.
\]
For $0<s<\|u\|_\infty$, the function $u+s$ is a continuous negative plurisubharmonic exhaustion of $E_s$, so $E_s$ is hyperconvex.
Integrating over the sublevel sets gives
\[
 (n+1)\mathcal E(u)=\int_B(-u)(\ddc u)^n
 =\int_0^{\|u\|_\infty}m(s)\,ds.
\]
For $0<s<t<\|u\|_\infty$, set $v=\frac{u+s}{t-s}$ and let $h_{s,t}$ be the relative extremal function of $F_t$ in $E_s$. Since $v\leq0$ on $E_s$ and $v\leq-1$ on $F_t$, the comparison principle gives
\begin{eqnarray*}
 (t-s)^n\Capa(F_t,E_s)
 &=&(t-s)^n\int_{E_s}(\ddc h_{s,t})^n\\
 &\leq&(t-s)^n\int_{E_s}(\ddc v)^n
 =\int_{E_s}(\ddc u)^n=m(s).
\end{eqnarray*}
Thus
\begin{equation}\label{eq:relative-sublevel-capacity}
 \Capa(F_t,E_s)^{-\frac{1}{n}}
 \geq\frac{t-s}{m(s)^{\frac{1}{n}}}.
\end{equation}
Combining \eqref{eq:relative-sublevel-capacity} with Lemma~\ref{lem:nested-capacity} gives
\[
 \Capa(F_t,B)^{-\frac{1}{n}}-\Capa(F_s,B)^{-\frac{1}{n}}
 \geq\frac{t-s}{m(s)^{\frac{1}{n}}}.
\]
Here we also used the monotonicity of capacity.
Set $r(s)=\Capa(F_s,B)^{-\frac{1}{n}}$. Then
\begin{equation}\label{eq:capacity-increment}
 r(t)-r(s)\geq\frac{t-s}{m(s)^{\frac{1}{n}}}.
\end{equation}
If $r$ is smooth, this gives $r'(s)\geq m(s)^{-\frac{1}{n}}$. For its inverse $w$, extended constantly beyond its range, the change of variables $y=r(t)$ yields
\[
 \int_0^\infty|w'(y)|^{n+1}\,dy
 =\int_0^{\|u\|_\infty}\frac{dt}{r'(t)^n}
 \leq\int_0^{\|u\|_\infty}m(t)\,dt
 =(n+1)\mathcal E(u).
\]
In general, define
\[
 w(y)=\sup\bigl(\{0\}\cup
 \{t\in(0,\|u\|_\infty):r(t)\leq y\}\bigr).
\]
Since $m(s)\leq\int_B(\ddc u)^n<\infty$,
\eqref{eq:capacity-increment} shows that $w(0)=0$, $w$ is Lipschitz, and
$0\leq w'\leq m(w)^{\frac{1}{n}}$ almost everywhere.
Moreover, $w(y)\to\|u\|_\infty$ as $y\to\infty$.
The change of variables for the nondecreasing absolutely continuous
function $w$ therefore gives
\begin{equation}\label{eq:inverse-capacity-energy}
 \int_0^\infty|w'|^{n+1}\,dy
 \leq\int_0^\infty m(w(y))w'(y)\,dy
 =\int_0^{\|u\|_\infty}m(t)\,dt
 =(n+1)\mathcal E(u).
\end{equation}
The volume--capacity inequality \eqref{eq:volume-capacity} gives
\[
 \mu\{-u\geq t\}
 \leq A\exp\left[-\beta\Capa(F_t,B)^{-\frac{1}{n}}\right]
 =Ae^{-\beta r(t)}.
\]
By the definition of $w$, $y\geq r(t)$ implies $w(y)\geq t$. Since
\[
 e^{-\beta r(t)}=\int_{r(t)}^\infty\beta e^{-\beta y}\,dy,
\]
we obtain
\[
 \mu\{-u\geq t\}\leq A\int_{\{w\geq t\}}\beta e^{-\beta y}\,dy.
\]
Also, \eqref{eq:volume-capacity} implies $\mu(B)\leq A$. Integrating the preceding estimate therefore gives, for every nonnegative increasing $C^1$ function $F$,
\[
 \int_BF(-u)\,d\mu
 \leq A\int_0^\infty F(w(y))\beta e^{-\beta y}\,dy.
\]
Take
\[
 F(t)=\exp\left[\beta\frac{t^{\frac{n+1}{n}}}
 {[(n+1)\mathcal E(u)]^{\frac{1}{n}}}\right].
\]
To apply Lemma~\ref{lem:Moser}, set
\[
 \widetilde w(t)
 =-\frac{\beta^{\frac{n}{n+1}}}
 {[(n+1)\mathcal E(u)]^{\frac{1}{n+1}}}
 w\left(-\frac{t}{\beta}\right),\qquad t\leq0.
\]
By \eqref{eq:inverse-capacity-energy}, $\displaystyle\int_{-\infty}^0|\widetilde w'|^{n+1}\,dt\leq1$. The lemma and the change of variables $t=-\beta y$ now give
\[
 \int_B\exp\left[\beta\frac{(-u)^{\frac{n+1}{n}}}
 {[(n+1)\mathcal E(u)]^{\frac{1}{n}}}\right]d\mu
 \leq AC_n.
\]
This proves the proposition.
\end{proof}
\begin{remark}
The same argument gives a weighted-energy version of
Proposition~\ref{prop:capacity-to-MT}. For $p>1$ and
$0\ne u\in\mathcal H_0(B)$, define
\[
 \mathcal E_p(u)=\frac{1}{n+1}\int_B(-u)^p(\ddc u)^n.
\]
With $w$ and $m$ as above, we have
\[
\begin{aligned}
p\int_0^\infty w(y)^{p-1}|w'(y)|^{n+1}\,dy
&\le p\int_0^\infty w(y)^{p-1}m(w(y))w'(y)\,dy\\
&=p\int_0^{\|u\|_\infty}s^{p-1}m(s)\,ds\\
&=\int_B(-u)^p(\ddc u)^n=(n+1)\mathcal E_p(u).
\end{aligned}
\]
Consequently,
\[
 \int_0^\infty
 \left|\left(w^{\frac{n+p}{n+1}}\right)'\right|^{n+1}\,dy
 \leq\frac{(n+p)^{n+1}}{p(n+1)^n}\mathcal E_p(u).
\]
Applying Lemma~\ref{lem:Moser} to the normalized function
$w^{\frac{n+p}{n+1}}$, with the same change of variables as above, gives
under assumption~\eqref{eq:volume-capacity}
\[
 \int_B\exp\left[
 \frac{\beta(n+1)p^{\frac{1}{n}}}{(n+p)^{\frac{n+1}{n}}}
 \frac{(-u)^{\frac{n+p}{n}}}{\mathcal E_p(u)^{\frac{1}{n}}}
 \right]d\mu\leq AC_n.
\]
In particular, for Euclidean volume on the ball,
\[
\int_B
\exp\!\left(
\alpha_{n,p}\frac{(-u)^{\frac{n+p}{n}}}{\mathcal E_p(u)^{\frac{1}{n}}}
\right)dV
\le C_{n,p},
\]
where
\[
 \alpha_{n,p}
 =\frac{2\pi n(n+1)p^{\frac{1}{n}}}{(n+p)^{\frac{n+1}{n}}}.
\]
For $p=1$, the same computation recovers the coefficient $\alpha_n$
and the energy $\mathcal E(u)$.

\end{remark}

\section{Moser--Trudinger inequalities on general domains}\label{sec:general-domains}

With the inequalities on the ball established, we now prove domain
monotonicity and complete the proof of Theorem~\ref{thm:sharp-MT}.
Comparison with an enclosing ball gives the critical exponential estimate.
The same argument transfers the logarithmic inequality to general domains.

\subsection{Domain monotonicity}\label{sec:domain-monotonicity}

We first prove domain monotonicity for integral functionals with fixed
Monge--Amp\`ere energy, in the spirit of \cite[Lemma~3.1]{TiW10}.

For a nonnegative convex function $F\in C^1([0,\infty))$, write
$f=F'$ and assume that $f(t)>0$ for $t>0$. Set
\[
 Y_F(\Omega)=\sup\left\{\int_\Omega F(-u)\dV:
 u\in\mathcal H_0(\Omega),\ \mathcal E(u)=1\right\}.
\]
The value $+\infty$ is allowed. Here and below, the energy is computed
on the domain of the function under consideration. The classes
$\mathcal E_0(\Omega)$ and $\mathcal E^1(\Omega)$ are defined as in
Section~\ref{subsec:sharp-inequalities}, with $B$ replaced by $\Omega$.

\begin{lemma}\label{lem:domain-monotonicity}
Let $\Omega_1\subset\Omega_2\subset\C^n$ be bounded strictly
pseudoconvex domains with smooth boundaries. Then
\begin{equation}\label{eq:domain-monotonicity}
 Y_F(\Omega_1)\le Y_F(\Omega_2).
\end{equation}
\end{lemma}
\begin{proof}
Fix $u\in\mathcal H_0(\Omega_1)$ with $\mathcal E(u)=1$.
Then $f(-u)$ is a nonnegative bounded density with positive mass.
Solve
\[
 \begin{cases}
 (\ddc w_1)^n=f(-u)\,dV&\text{in }\Omega_1,\\
 w_1=0&\text{on }\partial\Omega_1,
 \end{cases}
 \qquad
 \begin{cases}
 (\ddc w_2)^n=\mathbf1_{\Omega_1}f(-u)\,dV&\text{in }\Omega_2,\\
 w_2=0&\text{on }\partial\Omega_2.
 \end{cases}
\]
The Dirichlet existence theorem \cite{K98} gives
$w_i\in\PSH(\Omega_i)\cap C(\overline{\Omega_i})$.
Their boundedness and finite Monge--Amp\`ere mass imply
$w_i\in\mathcal E_0(\Omega_i)$.
On $\Omega_1$ their Monge--Amp\`ere measures agree, and
$w_2\le0=w_1$ on $\partial\Omega_1$. The comparison principle
therefore gives $w_2\le w_1$ on $\Omega_1$.

Set $a_i=\mathcal E(w_i)>0$. Since the density of $w_2$ vanishes
outside $\Omega_1$,
\[
 (n+1)a_2=\int_{\Omega_1}(-w_2)f(-u)\dV
 \ge\int_{\Omega_1}(-w_1)f(-u)\dV=(n+1)a_1.
\]
Let $v=a_2^{-\frac{1}{n+1}}w_2$, so that $\mathcal E(v)=1$.
The mixed-energy inequality
\cite[Theorem~2.1, with $p=1$]{ACC12} reads
\begin{equation}\label{eq:domain-mixed-energy}
 \int_\Omega(-\psi)(\ddc w)^n
 \le(n+1)\mathcal E(\psi)^{\frac{1}{n+1}}
              \mathcal E(w)^{\frac{n}{n+1}},
 \qquad \psi,w\in\mathcal E^1(\Omega).
\end{equation}
For $n=1$ this is the Cauchy--Schwarz inequality for the Dirichlet
energy. Applying \eqref{eq:domain-mixed-energy} to $u,w_1$, we obtain
\[
 \int_{\Omega_1}(-u)f(-u)\dV
 \le(n+1)a_1^{\frac{n}{n+1}}
 \le(n+1)a_2^{\frac{n}{n+1}}
 =\int_{\Omega_1}(-v)f(-u)\dV.
\]
Consequently, 
\[\int_{\Omega_1}f(-u)(u-v)\dV\ge0.\]
Convexity gives
$F(-v)\ge F(-u)+f(-u)(u-v)$ on $\Omega_1$, and hence
\begin{equation}\label{eq:domain-functional-comparison}
 \int_{\Omega_2}F(-v)\dV
 \ge\int_{\Omega_1}F(-v)\dV
 \ge\int_{\Omega_1}F(-u)\dV.
\end{equation}

To obtain smooth competitors, choose positive functions
$f_j\in C^\infty(\overline{\Omega_2})$ converging to
$\mathbf1_{\Omega_1}f(-u)$ in $L^2(\Omega_2)$, with uniformly bounded
$L^\infty$ norms. The smooth Dirichlet solutions
\[
 (\ddc w_{2,j})^n=f_j\,dV,\qquad
 w_{2,j}|_{\partial\Omega_2}=0
\]
exist by \cite{CKNS85} and converge uniformly to $w_2$ by stability;
see also the approximation in \cite[Lemma~3.1]{WWZ20}.
It follows that
\[
 \mathcal E(w_{2,j})
 =\frac{1}{n+1}\int_{\Omega_2}(-w_{2,j})f_j\dV\longrightarrow a_2.
\]
Thus $v_j=\mathcal E(w_{2,j})^{-\frac{1}{n+1}}w_{2,j}$ belongs to
$\mathcal H_0(\Omega_2)$, has energy one, and converges uniformly
to $v$. Passing to the limit in $\displaystyle\int_{\Omega_2}F(-v_j)\dV$ and
using \eqref{eq:domain-functional-comparison}, we get
\[
 Y_F(\Omega_2)\ge\int_{\Omega_1}F(-u)\dV.
\]
Taking the supremum over $u$ proves \eqref{eq:domain-monotonicity}.
\end{proof}

\begin{corollary}\label{cor:mt-ball-transfer}
For $\alpha>0$ and $\beta\ge1$, define
\[
 Y_{\alpha,\beta}(\Omega)
 =\sup_{0\ne u\in\mathcal H_0(\Omega)}
 \int_\Omega\exp\left[\alpha
 \left(\frac{-u}{\mathcal E(u)^{\frac{1}{n+1}}}\right)^\beta\right]\dV.
\]
If $B_r(z_0)\subset\Omega\subset B_R(z_1)$, then
\begin{equation}\label{eq:mt-domain-sandwich}
 r^{2n}Y_{\alpha,\beta}(B)
 \le Y_{\alpha,\beta}(\Omega)
 \le R^{2n}Y_{\alpha,\beta}(B).
\end{equation}
In particular, $Y_{\alpha,\beta}(\Omega)$ is finite if and only if $Y_{\alpha,\beta}(B)$ is finite, also for $\alpha=\alpha_n$ and $\beta=\frac{n+1}{n}$.
\end{corollary}
\begin{proof}
Apply Lemma~\ref{lem:domain-monotonicity} with
$F(t)=e^{\alpha t^\beta}$. For $u$ on $B_s(z_0)$, the change of
variables $z=z_0+s\zeta$ gives
\[
 \mathcal E\bigl(u(z_0+s\,\cdot)\bigr)=\mathcal E(u),
 \qquad
 Y_{\alpha,\beta}(B_s(z_0))=s^{2n}Y_{\alpha,\beta}(B).
\]
Domain monotonicity applied to the two inclusions proves
\eqref{eq:mt-domain-sandwich}.
\end{proof}

\subsection{The inequalities on general domains}

\begin{corollary}\label{cor:general-log}
Let $\Omega\Subset B_R(a)\subset\mathbb C^n$ be a strictly pseudoconvex
domain with smooth boundary. Then, for every $u\in\mathcal H_0(\Omega)$,
\[
 \log\left(\frac{1}{|\Omega|}\int_\Omega e^{-u}\dV\right)
 \leq\frac{\mathcal E(u)}{[2\pi(n+1)]^n}
 +\sum_{j=1}^n\frac{1}{j}+\log\frac{|B_R|}{|\Omega|}.
\]
The energy coefficient is optimal, allowing an arbitrary finite additive
constant depending on $\Omega$.
\end{corollary}
\begin{proof}
For $t>0$, put $F_t(s)=e^{ts}$. By Theorem~\ref{thm:sharp-log},
\[
 Y_{F_t}(B)\leq|B|\cdot\exp\left(
 \frac{t^{n+1}}{[2\pi(n+1)]^n}+\sum_{j=1}^n\frac{1}{j}\right).
\]
Lemma~\ref{lem:domain-monotonicity} and spatial scaling give
$Y_{F_t}(\Omega)\leq R^{2n}Y_{F_t}(B)$.
For $u\ne0$, take $t=\mathcal E(u)^{\frac{1}{n+1}}$ and apply this bound
to $\frac{u}{t}$, which has energy one. This proves the inequality;
the case $u=0$ is immediate.

For optimality, choose $B_r(z_0)\subset\Omega$. If a logarithmic inequality
held on $\Omega$ with a smaller energy coefficient $A$ and a finite
additive constant $C_\Omega$, then
\[
 r^{2n}Y_{F_t}(B)\leq Y_{F_t}(\Omega)
 \leq|\Omega|\exp(At^{n+1}+C_\Omega),\qquad t>0.
\]
Rescaling each nonzero $v\in\mathcal H_0(B)$ to have energy one would give
a logarithmic inequality on $B$ with the same smaller energy coefficient.
This contradicts Theorem~\ref{thm:sharp-log}.
\end{proof}

\begin{proof}[Proof of Theorem~\ref{thm:sharp-MT}]
Apply Proposition~\ref{prop:capacity-to-MT} with
$d\mu=\frac{dV}{|B|}$,
$A=\exp\left(\sum_{j=1}^n\frac{1}{j}\right)$ and $\beta=2\pi n$.
Since $\frac{2\pi n}{(n+1)^{\frac{1}{n}}}=\alpha_n$, we obtain
\[
 \int_B\exp\left[\alpha_n
 \left(\frac{-u}{\mathcal E(u)^{\frac{1}{n+1}}}\right)^{\frac{n+1}{n}}\right]\dV
 \leq C_n|B|\cdot\exp\left(\sum_{j=1}^n\frac{1}{j}\right).
\]
Corollary~\ref{cor:mt-ball-transfer} now gives the stated estimate on $\Omega$.
For optimality, an exponential inequality with $\alpha>\alpha_n$ would imply
\eqref{eq:log-energy-mt} with energy coefficient
\[
 \frac{n^n}{(n+1)^{n+1}\alpha^n}<\frac{1}{[2\pi(n+1)]^n}.
\]
This contradicts Corollary~\ref{cor:general-log}.
\end{proof}

\begin{remark}
For comparison, tracking the parameter dependence in the quasi-sharp
inequalities of \cite[Sections~4--5]{BB22} gives
\[
 \int_\Omega\exp\left[(1-\varepsilon)\alpha_n
 \left(\frac{-u}{\|u\|_{\mathcal H_0(\Omega)}}\right)^{\frac{n+1}{n}}\right]\dV
 \leq\frac{C(n,\Omega)}{\varepsilon^n},\qquad 0<\varepsilon<1.
\]
Theorem~\ref{thm:sharp-MT} gives a bound independent of $\varepsilon$
and includes the value $\varepsilon=0$.
\end{remark}

\end{document}